\documentclass[11pt,a4paper,final]{amsart} 
\usepackage[margin=24mm]{geometry}
\usepackage{anyfontsize}
\usepackage{colonequals}
\usepackage{mathtools}
\usepackage{amsmath}
\usepackage{amssymb}
\usepackage{multicol}

\usepackage[switch,mathlines]{lineno}

\usepackage[skip=1em,font=footnotesize]{caption}
\usepackage{subcaption}

\usepackage[inline]{enumitem}
\usepackage[disableredefinitions]{complexity}
\usepackage[usenames,dvipsnames]{xcolor}
\usepackage[pdftex,colorlinks,citecolor=blue,urlcolor=blue]{hyperref}
\usepackage[nameinlink]{cleveref}
\usepackage{array}
\usepackage{tabularx}

\usepackage{comment}
\usepackage[author={AD. Santamaria},icon=Comment,color=SpringGreen!50,opacity=0.7]{pdfcomment}

\usepackage{graphicx}
\usepackage{wrapfig}

\numberwithin{equation}{section}
\newtheorem{theorem}{Theorem}
\newtheorem{proposition}[theorem]{Proposition}

\newtheorem{corollary}[theorem]{Corollary}

\newtheorem{lemma}[theorem]{Lemma}

\theoremstyle{definition}
\newtheorem{definition}[theorem]{Definition}
\newtheorem{example}[theorem]{Example}
\newtheorem{remark}[theorem]{Remark}

\newcommand{\introthmname}{}
\newtheorem{introthminn}{\introthmname}

\DeclareMathOperator{\lk}{\mathrm{link}}

\DeclareMathOperator{\del}{\mathrm{del}}
\DeclareMathOperator{\cl}{\mathrm{cl}}

\newcommand{\RR}{{\mathbb R}}

	\newcommand{\SpaceX}			{X}
	\newcommand{\Disk}			{D}
	\newcommand{\Cycle}			{S^1}
	\newcommand{\Mobius}			{M}
	\newcommand{\Annulus}		{A}
	\newcommand{\Ball}			{B}
	\newcommand{\Pplane}			{\RR\mathrm{P}^2}
	\newcommand{\Pspace}[1]			{\RR\mathrm{P}^{#1}}

	\newcommand{\DunceH}			{Z}
	
	\newcommand{\dMobius}		{\partial\Mobius}
	
	\newcommand{\Triang}[1]		{\Delta_{#1}}
	\newcommand{\TriangD}		{\Triang{\Disk}}
	\newcommand{\TriangDi}[1]	{\Triang{\Disk_{#1}}}
	\newcommand{\TriangMo}		{\Triang{\Mobius}}
	\newcommand{\TriangAn}		{\Triang{\Annulus}}

	\newcommand{\TriangDH}		{\Triang{\DunceH}}
	\newcommand{\TriangPp}		{\Triang{\Pplane}}
	\newcommand{\TriangX}		{\Triang{\SpaceX}}
	
	\newcommand{\dTriangMo}		{\Triang{\dMobius}}

	\newcommand{\gSet}[1]		{\left\langle{#1}\right\rangle}

\title[Partitioning $\Pplane$ and $\DunceH$]
		{Partitioning the projective plane and the dunce hat}

\author[A.D. Santamaría-Galvis]{Andrés D. Santamaría-Galvis}
\thanks{The work of the author is supported in part by the Slovenian Research Agency project number N1-0160}
\address{
	Univerza na Primorskem,
	UP FAMNIT,
	Glagoljaška 8,
	Koper 6000, Slovenija
	}
\email{andres.santamaria@famnit.upr.si}
\urladdr{\url{https://sites.google.com/view/adsantamaria}}

\date{\today}


\begin{document}


	\begin{abstract}
		\label{sec:abstract}


		We show that any finite triangulation of the real projective plane or the dunce hat is partitionable.
		To prove this, we introduce simple yet useful gluing tools that allow us to reduce partitionability of a given complex to that of smaller constituent relative subcomplexes such as the disk or the open Möbius strip.
		The gluing process generates partitioning schemes with a distinctive shelling/constructible flavor. 
		We also give a tool to lift partitionability of relative simplicial complexes to that of the preimages of certain simplicial quotient maps.
		
	\end{abstract}
	\maketitle

	\section{Introduction}
		\label{sec:introduction}

		Partitionability (see \Cref{sec:preliminaries} for further details) is a combinatorial property of simplicial complexes that arose in the 1970s \cite{Ball:1977} \cite[App. 4]{Provan:1977} and appears early in the study of shellable and Cohen-Macaulay complexes \cite{Bjorner/Garsia/Stanley:1982,Stanley:1979}.
		It is known from its occurrence in the Stanley Partitioning Conjecture, which asks whether every Cohen-Macaulay complex is partitionable \cite[p. 149]{Stanley:1979}.
		Garsia asks for a similar question for order complexes in \cite[Remark 5.2]{Garsia:1980}. 
		If this conjecture were true, it would have explained combinatorially why the $h$-vector has non-negative entries.
		The fact that shellable simplicial complexes are both Cohen-Macaulay and partitionable seems to reinforce the idea.
		However, the Cohen-Macaulay yet non-partitionable counterexample constructed by Duval, Goeckner, Klivans and Martin in \cite{Duval/Goeckner/Klivans/Martin:2016} (see also \cite{Duval/Klivans/Martin:2017}) disproves the conjecture.

		Although \cite{Duval/Goeckner/Klivans/Martin:2016} settles the Stanley Partitionability Conjecture, it raises some other questions.
		For instance, we still do not know if the Stanley Conjecture holds when restricted to $2$-dimensional complexes.
		This question motivates the current paper, where we show that all triangulations of certain well known $2$-dimensional, non-shellable, Cohen-Macaulay spaces are partitionable.

		Our techniques are complementary to that of \cite{Duval/Goeckner/Klivans/Martin:2016}, where the authors built their counterexamples by gluing together several shellable $3$-balls along a common shellable $2$-ball.
		Their gluing process carries over the Cohen-Macaulay property from the original complexes to the new one, but it breaks partitionability.
		Their proof uses the framework of relative simplicial complexes to put restrictions on a partitioning scheme for the counterexample complex, eventually arriving at a contradiction.

		In contrast, we will develop tools in \Cref{stt:gluing_lemma,stt:assymetrical_gluing_lemma,stt:folding_lemma} that give sufficient conditions for producing partitioning schemes on and from constituent relative complexes.
		The first two lemmas carry over partitioning schemes of two relative complexes to a larger one, while the last lemma gives conditions for a quotient map to preserve a partitioning scheme.
		These lemmas will allow us to use a ``splitting-then-gluing-back'' strategy that improves the naïve facet-by-facet approach.
		Indeed, these techniques can be seen as a broad generalization of the partitioning aspects of shellability.

		We use the splitting and gluing approach on triangulations of the real projective plane $\Pplane$ and of the dunce hat $\DunceH$.
		Although particular triangulations of these spaces have long been known to be partitionable (see \cite{Hachimori:web}), we prove a considerably more general result.

		\begin{theorem}
			If $\Delta$ is any triangulation of the real projective plane $\Pplane$ or of the dunce hat $\DunceH$, then $\Delta$ is partitionable.
			\label{stt:main_result}
		\end{theorem}

		As an intermediate result, and using the same machinery, we prove partitionability of any triangulation of the open Möbius strip.
		Here the word \emph{open} means that we work relative to the boundary of the space.

		This paper is structured as follows.
		In \Cref{sec:preliminaries} we give a basic background on partitionability, relative complexes and related notions.
		The main tools we use are introduced and developed in \Cref{sec:partitioning_preserving_tools}.
		In \Cref{sec:partitionable_surfaces} we show that the $2$-disk and the Möbius strip are partitionable relative to certain portions of their respective boundaries.
		Finally, in \Cref{sec:main_theorem} we prove \Cref{stt:main_result}.


	\numberwithin{theorem}{section}

	\section{Preliminaries}
		\label{sec:preliminaries}

		In this section we succinctly present the main notions we require throughout the paper.
		We refer to \cite{Bjorner:1995, Stanley:1996} for additional background and notational conventions on simplicial combinatorics.

		An \emph{(abstract) simplicial complex} $\Delta$ is a collection of subsets closed under inclusion whose elements are taken from a finite ground set.
		For its standard geometric realization we use the notation $|\Delta|$.
		The elements of $\Delta$ are called \emph{faces} and the maximal faces are called \emph{facets}.
		By abuse of notation, when no confusion is induced, we denote the face $\{v_{1},v_{2},\dotsc,v_{k}\}$ as $v_{1}v_{2} \dotso v_{k}$, for $k \geq 1$.
		To make this more visible in our concrete examples, we will use monospaced typeface numerals.
		As usual, faces with one, two and three elements are called \emph{vertices}, \emph{edges} and \emph{triangles}, respectively.
		The collection \emph{generated} by a family of sets $\Lambda$ consists of all subsets of the sets in $\Lambda$, and it is denoted by $\gSet{\Lambda}$.
		This collection is a simplicial complex, sometimes called the \emph{combinatorial closure} of $\Lambda$.

		The \emph{dimension} $\dim \sigma$ of a face $\sigma \in \Delta$ is one less than the cardinality of $\sigma$, and the \emph{dimension} of $\Delta$, denoted as $\dim \Delta$, is the maximum dimension among the faces of $\Delta$.
		A complex is called \emph{pure} if all its facets have the same dimension.
		We are mostly interested in pure complexes, although some of the results in this paper will also hold for non-pure complexes.

		The \emph{link} of a face $\sigma \in \Delta$ is the complex defined by
		$\lk_\Delta(\sigma)
			\colonequals
		\{\tau \in \Delta \, : \, \tau \cap \sigma = \varnothing,\, \tau \cup \sigma \in \Delta \}$,
		and the \emph{deletion} of $\sigma$ is the complex
		$\del_\Delta(\sigma)
			\colonequals
		\{\tau \in \Delta \, : \, \tau \cap \sigma = \varnothing \}$.
		The \emph{cone} of a complex with \emph{apex} $u$ is defined by
		$u*\Delta
			\colonequals
		\{\{u\} \cup \sigma \, : \,  \sigma \in \Delta\} \, \cup \, \Delta$, $\{u\} \notin \Delta$.

		\subsection{Partitionability and shellability}
			\label{sec:partitionability_and_shellability}

			A simplicial complex admits a partial order by inclusion in a natural way, and we regard a simplicial complex as equivalent to its \emph{face poset}.
			We are interested in the problem of partitioning $\Delta$ (i.e.\ its face poset) into intervals of the form
			$\left[R(\sigma),\sigma\right]
				\colonequals
			\{\rho \in \Delta \, : \, R(\sigma) \subseteq \rho \subseteq \sigma,\; \text{where } \sigma \text{ is a facet in } \Delta \}$.
			Notice that an interval in the face poset of a simplicial complex is isomorphic to the Boolean lattice $B_{n}$, where $n$, the rank of $B_{n}$, is the cardinality of $\sigma \setminus R(\sigma)$.
			If $\Delta$ can be partitioned in such a way, we say that $\Delta$ is \emph{partitionable}, and the set of intervals is a \emph{partitioning scheme} (or a \emph{partitioning}, for short) of $\Delta$.
			\Cref{fig:graph_example} shows a partitionable simplicial complex with a partitioning scheme.

			The \emph{$f$-vector} of a $(d-1)$-dimensional complex $\Delta$ is
			$f(\Delta) \colonequals \left(f_{-1}, f_{0}, \dotsc, f_{d-1} \right)$,
			where $f_i$ denotes the number of $i$-dimensional faces of $\Delta$.
			The entry $f_{-1}$ is $0$ if $\Delta$ is empty and $1$ otherwise.
			The \emph{(reduced) Euler characteristic} of $\Delta$ is defined as
			$\tilde{\chi}(\Delta) \colonequals \sum\nolimits_{i=0}^{d} {(-1)^{i-1} f_{i-1}}$.
			This value is a topological invariant of the space that $\Delta$ triangulates, so the choice of triangulation does not modify the Euler characteristic.
			The $h$-vector of $\Delta$ is
			$h(\Delta) \colonequals \left(h_{0}, h_{1}, \dotsc, h_{d}\right)$
			where its entries are obtained from the relation
			\begin{equation}
				\sum_{i=0}^{d} {f_{i-1} (t-1)^{d-i}} = \sum_{i=0}^{d} {h_{i} t^{d-i}}.
				\label{equ:h_pol_and_f_pol}
			\end{equation}
			In particular, $h_d = (-1)^{d-1}\tilde{\chi} (\Delta)$ and $f_{d-1} = \sum\nolimits_{i=0}^{d} h_i$.
			The $h$-vector of a pure partitionable complex has a combinatorial meaning: the $i$th entry of $h(\Delta)$ counts how many intervals in a partitioning of $\Delta$ have a minimal face $R(\cdot)$ with dimension $i-1$.
			This is true regardless of the chosen partitioning scheme.

			\begin{example}
				The graph ($1$-dimensional complex) appearing in \Cref{fig:graph_example} (a) has $h$-vector $\left(1, 3, 0\right)$.
				Thus, if the graph is partitionable, any partitioning must match one edge to the empty set, three edges to one incident vertex each, and no edge matches to itself.
				In this case, a partitioning exists and it is shown in \Cref{fig:graph_example} (b) and (c).

				\begin{figure}[!ht]
					\begin{center}
						\centering
						\begin{tabular}{	>{\centering\arraybackslash}m{4.5cm}
												>{\centering\arraybackslash}m{5.0cm}
												>{\centering\arraybackslash}m{3.0cm}
												}
							\vspace{ 0.20em} \includegraphics[scale=0.40]{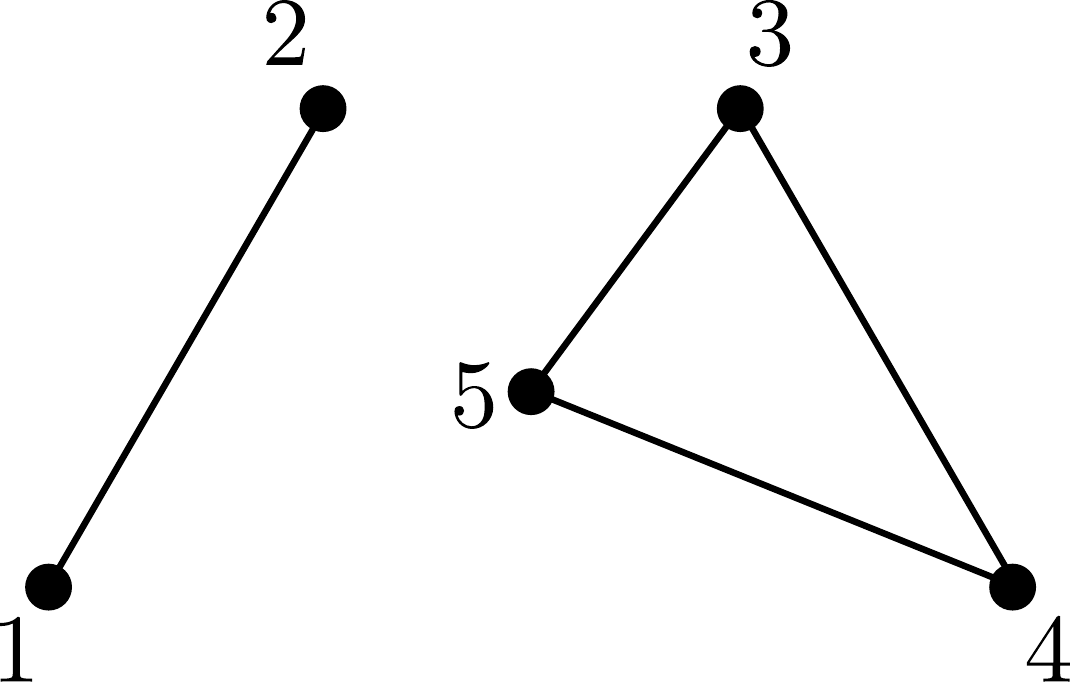} &
							\vspace{ 0.00em} \includegraphics[scale=0.46]{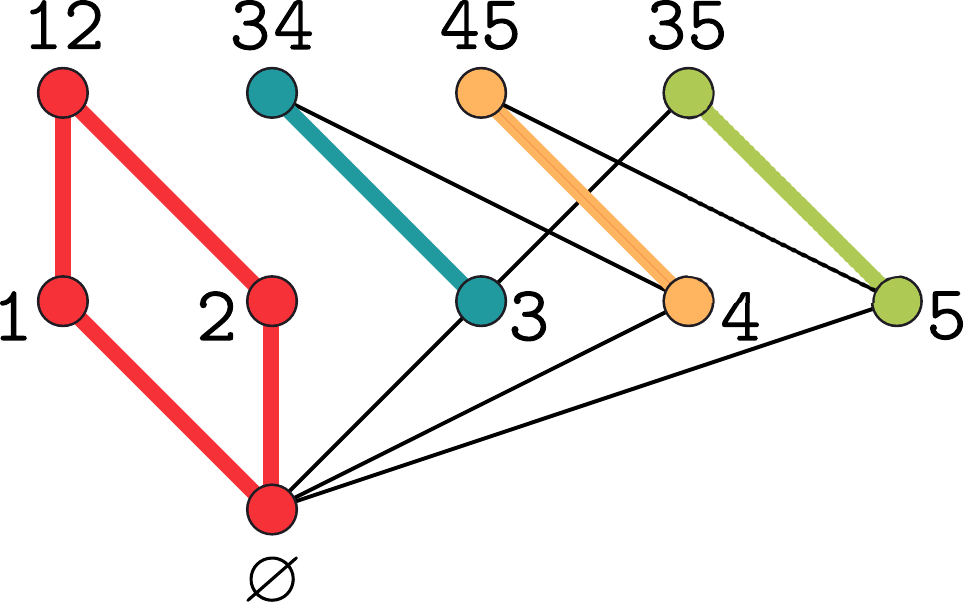} &
							{$\!\begin{aligned}
								\left[ \mathtt{\varnothing, 12} \right] \ & \sqcup \\
								\left[           \mathtt{3, 34} \right] \ & \sqcup \\
								\left[           \mathtt{4, 45} \right] \ & \sqcup \\
								\left[           \mathtt{5, 35} \right] \ &
							\end{aligned}$}
							\vspace{ 0.20em}
							\\
							(a) & (b) & (c)
						\end{tabular}
						\caption{
							A partitionable $1$-dimensional complex (a).
							Its partitioning scheme is depicted in (b) and expressed as the union of intervals in (c).}
						\label{fig:graph_example}
					\end{center}
				\end{figure}
				In general, it is known that a graph is partitionable if and only if at most one of its connected components is a tree \cite[Problem 35]{Kaibel/Pfetsch:2003}.
				See also \Cref{stt:part_graphs_resp_empty_face}.
				\label{stt:graphs_are_partitionable_ex}
			\end{example}

			\begin{example}
				The \emph{bow-tie} complex $\Delta = \gSet{\mathtt{125, 345}}$ is non-partitionable since $h(\Delta) = (1,2,-1,0)$ has negative entries.
				\label{stt:bow-tie_partitionable_ex}
			\end{example}

			A notion related to partitionability is that of shellability.
			A simplicial complex $\Delta$ is \emph{shellable} if its facets can be ordered linearly in such a way that the complex generated by the facet $\sigma_{j}$, for each $j > 1$, intersects in a pure $(\dim\sigma_j-1)$-dimensional complex with the complex generated by the previous $j-1$ facets in the ordering.
			The $j$th step in the shelling process adds a unique minimal new face $R(\sigma_{j})$.
			Thus, every time a facet $\sigma_{j}$ is attached to the shelling it introduces the interval $[R(\sigma_{j}), \sigma_{j}]$ with just new faces.
			Hence, the shelling process inductively produces a partitioning scheme of $\Delta$.
			In a recent work, Hachimori defines several properties satisfied by partitioning schemes induced from shellings \cite{Hachimori:2021}.
			
			Although every shellable simplicial complex is partitionable, the converse is not true.
			\Cref{stt:graphs_are_partitionable_ex} shows a partitionable but non-shellable simplicial complex.
			In general, it is known that a graph is shellable if and only if it is connected \cite[Problem 34]{Kaibel/Pfetsch:2003}.
			See a related result in \Cref{stt:shell_graphs_resp_empty_face}.
					
			It is also useful to remember that the cone operation preserves shellability and partitionability.
			Furthermore, a complex is partitionable (shellable) if and only if its cone is.

		\subsection{Relative complexes}
			\label{sec:relative_complexes}

			A \emph{relative simplicial complex} $\Phi = \left(\Delta,\Gamma\right)$ is a pair of simplicial complexes, where $\Gamma$ is a subcomplex of $\Delta$.
			If $\Gamma \neq \varnothing$, the relative complex $\Phi$ is said to be \emph{proper}.
			Lacking a better way to refer to $\Gamma$, we say that $\Gamma$ is the \emph{relative part} of $\Phi$.
			We regard the usual (non-relative) simplicial complex as the relative simplicial complex with $\Gamma = \varnothing$.
			The \emph{faces} of $\left(\Delta,\Gamma\right)$ are the faces of $\Delta$ that are not in $\Gamma$, and the maximal faces are the \emph{facets} of $\left(\Delta,\Gamma\right)$.
			We say that $\left(\Delta,\Gamma\right)$ is \emph{pure} if its facets have the same dimension.
			We use the frequently overloaded symbol $\Delta \setminus \Gamma$ to denote the set of faces of $\left(\Delta,\Gamma\right)$, namely by taking $\Delta$ and $\Gamma$ as set systems.
			We will use this notation even if $\Gamma$ is not a subcomplex of $\Delta$.
			We define the \emph{minimal representation} of $\left(\Delta,\Gamma\right)$ as the complex $\left(\gSet{\Lambda},\gSet{\Lambda} \setminus \Lambda \right)$
			where $\Lambda \colonequals \Delta \setminus \Gamma$.

			The relative complex $\left( \Delta,\Gamma \right)$ can be seen as a combinatorial model of the topological pair of spaces $\left(|\Delta| , |\Gamma| \right)$
			(see, e.g., \cite[\S 2.1--2.2]{Hatcher:2002} and \cite[\S 1.9]{Munkres:1984}).
			This relative notion is used in algebraic topology to calculate the homology of the quotient space $|\Delta| \, / \, |\Gamma|$.
			Furthermore, relative complexes offer a natural setting to build up new complexes with a given property, as in \cite{Duval:1999,Duval:2000,Duval/Goeckner/Klivans/Martin:2016,Juhnke-Kubitzke/Venturello:2019,Santamaria/Woodroofe:2021}.

			The combinatorics of non-relative simplicial complexes carries over straightforwardly to their relative counterparts, including notions such as face poset and partitionability.
			The \emph{face poset} of the relative complex $\Phi = \left( \Delta, \Gamma \right)$ is the set system $\Delta \setminus \Gamma$ ordered by inclusion.
			This partial order of the elements of $\Phi$ with respect to the inclusion is order-convex, i.e.\ if $\sigma, \tau \in \Phi$, and $\sigma \subseteq \rho \subseteq \tau$, then $\rho \in \Phi$.
			A relative complex $\Phi = \left( \Delta, \Gamma \right)$ is said to be \emph{partitionable} if its face poset can be written as the disjoint union of intervals in such a way that the top element of each is a facet of $\Phi$.
			We also say that $\Delta$ \emph{is partitionable relative to} (or \emph{with respect to}) $\Gamma$.
			Such a disjoint union of intervals is a \emph{partitioning scheme} of the relative complex (or a \emph{partitioning}, for short).
			\Cref{fig:isomorphic_posets} shows a partitionable relative simplicial complex (d) with a partitioning scheme drawn over its face poset (e) and written as a disjoint union of intervals (f).
			Similar to the non-relative setting, the \emph{$f$-vector} of a relative complex $f(\Delta,\Gamma)$ stores the number of faces of a given dimension of $\left(\Delta,\Gamma\right)$, and the \emph{$h$-vector} $h(\Delta,\Gamma)$ is determined by the $f$-vector through \labelcref{equ:h_pol_and_f_pol}.
			When $(\Delta,\Gamma)$ is pure, the $h$-vector has the same combinatorial meaning as its non-relative pure counterpart, namely the $i$th entry counts the number of intervals in a partitioning scheme of $\left(\Delta,\Gamma\right)$ having minimal face $R(\cdot)$ with dimension $i - 1$.

			Shellability can also be defined for relative complexes.
			If $\Gamma = \varnothing$, we have the usual non-relative notion.
			Let $\Gamma \neq \varnothing$.
			An ordering $\sigma_{1}, \dotsc, \sigma_{m}$ of the facets of the complex $\left(\Delta,\Gamma\right)$ is a \emph{shelling} if, for each $j \geq 1$, the facet $\sigma_{j}$ intersects in a pure $(\dim\sigma_j-1)$-dimensional complex with the complex $\Lambda_{j-1}$
			generated by $\sigma_{1}, \dotsc, \sigma_{j-1}$ together with $\Gamma$. 
			Set $\Lambda_{0} \colonequals \Gamma$.
			If a shelling exists, we say that $\left(\Delta,\Gamma\right)$ is \emph{shellable}, or that $\Delta$ \emph{is shellable relative to} (or \emph{with respect to}) $\Gamma$.
			From the definition, it is clear that the $j$th step of the shelling attaches all the faces from
			$\langle \sigma_{j} \rangle \setminus \Lambda_{j-1}$ to $\Lambda_{j-1}$, and a unique minimal face $R(\sigma_{j}) \in \sigma_{j}$ is included into the shelling (just as in its non-relative counterpart).
			Thus $\langle \sigma_{j} \rangle \setminus \Lambda_{j-1} = [R(\sigma_{j}), \sigma_{j}]$.
			As a consequence, a shellable relative simplicial complex is also partitionable.

			\begin{remark}
				We digress here with an algorithmic description for the construction of partitioning schemes of relative simplicial complexes that fits with the approach of this paper.

				Consider a relative simplicial complex $\left(\Delta, \Gamma \right)$.
				We could think of a partitioning scheme as the outcome of a procedure that constantly takes facets from $\Delta$ to place them into $\Gamma$:

				Start a partitioning of $\left(\Delta, \Gamma \right)$ by choosing a suitable facet $\sigma$.
				Since continuing the partitioning is equivalent to considering the problem anew for the complex $\Phi = \left( \Delta, \Gamma \cup \gSet{\sigma} \right)$, we can naturally continue the process by recursively placing a new suitable facet into the relative part of $\Phi$.
				The process ends when the first part and the second part of the pair coincide.
				If any failure appears, backtrack and select the next suitable facet.
				See the proofs of \Cref{stt:rel_disk_is_partitionable} and \Cref{stt:links_of_balls_as_balls} for a related approach.

				Alternatively, we may consider the minimal representation of the complex every time a new facet is included.
				In this case, the procedure ends when we reach the empty complex.
				This alternative is also possible given the poset isomorphism between a relative simplicial complex and its minimal representation.
				See \Cref{stt:isoposets_partitionable} and \Cref{stt:conseq_isoposets_part_ex} (\labelcref{stt:part_minimal_representation}).
				
				We consider worth mentioning in this side note that an integer programming routine to decide partitionability has been implemented for SageMath \cite{sagemath}.
				\label{stt:algorithm_feeding_rel_part}
			\end{remark}

		\subsection{Triangulations on surfaces}
			\label{sec:triangulations_on_surfaces}

			In this work we are concerned with partitionability of finite triangulations of concrete topological spaces.
			We say that a \emph{(simplicial) triangulation} of a space $\SpaceX$ is a simplicial complex
			$\TriangX$ whose underlying space is $\SpaceX$, i.e.\ $|\TriangX| \cong X$.
			We relax the language to say that a topological space is \emph{partitionable} (or \emph{shellable}) when all of its finite triangulations are.
			Again, we extend these notions to their relative counterparts when required.

			In \Cref{sec:partitioning_pplane} we are going to need familiar results on surfaces with and without boundary.
			Given that every surface can be triangulated (see e.g.\ \cite{Thomassen:1992}), the homeomorphism between surfaces can be reduced to the combinatorial comparison of certain topological invariants provided by the classification theorem for compact surfaces.
			We refer the reader to \cite{Gallier/Xu:2013} for further details.


	\section{Tools that preserve partitioning schemes}
		\label{sec:partitioning_preserving_tools}

		We start this section with two lemmas that give conditions under which we can glue complexes while preserving partitioning schemes.
		They build a (relative) partitionable simplicial complex out of two partitionable relative simplicial complexes.

		\begin{lemma}[Gluing Lemma]
			Let $\Phi_{a} =  \left(\Delta_{a}, \Gamma_{a} \right)$ and
				 $\Phi_{b} =  \left(\Delta_{b}, \Gamma_{b} \right)$
			be partitionable complexes such that no facet of one is properly contained in a facet of the other, and let
			$\Sigma \colonequals
				\left( \Gamma_{a} \setminus \Delta_{b} \right)  \, \cup \,
				\left( \Gamma_{b} \setminus \Delta_{a} \right)  \, \cup \,
				\left( \Gamma_{a} \cap \Gamma_{b} \right)$.
			The complex
			$\Phi =  \left( \Delta_{a} \cup \Delta_{b}, \Sigma  \right)$
			is partitionable if the following conditions are met:
			\begin{enumerate}[label=\tt(G\arabic*)]
				\item \label{con:intersection_embedded}
					$\Gamma_{a} \cup \Gamma_{b} \supseteq \Delta_{a} \cap \Delta_{b}$.
				\item \label{con:partition}
					$\Sigma$ is a subcomplex of $\Delta_{a} \cup \Delta_{b}$.
					In particular, the set system $\Sigma$ is a simplicial complex.
			\end{enumerate}
			\label{stt:gluing_lemma}
		\end{lemma}

		\begin{proof}
			It is sufficient to show that the union of the partitioning schemes of $\Phi_{a}$ and $\Phi_{b}$ forms a partitioning scheme of $\Phi$.
			First, we want to prove that the faces of $\Phi_{a}$ are disjoint from those of $\Phi_{b}$.
			By \labelcref{con:intersection_embedded} and properties of set difference
			$\left(\Delta_{a} \setminus \Gamma_{a}\right) \, \cap \, \left(\Delta_{b} \setminus \Gamma_{b} \right) =
			 \left( \Delta_{a} \cap \Delta_{b} \right)  \, \setminus \,
			 \left( \Gamma_{a} \cup \Gamma_{b} \right) = \varnothing$.
			Thus, it is clear that no interval in the partitioning of $\Phi_{a}$ shares a face with an interval from $\Phi_{b}$.
			Moreover, since the facets of $\Phi_{a}$ and $\Phi_{b}$ are facets in $\Phi$, the partitioning schemes of each can be directly carried over into $\Phi$.

			Now we need to check that $\Phi$ is the union of the $\Phi_{a}$ and $\Phi_{b}$.
			Notice that any face in
			$\Lambda \colonequals \left(\Delta_{a} \setminus \Gamma_{a}\right) \, \cup \, \left(\Delta_{b} \setminus \Gamma_{b} \right)$
			is a face of $\Delta_{a} \cup \Delta_{b}$.
			Then, by elementary set theory, the faces of $\Delta_{a} \cup \Delta_{b}$ we are missing in $\Lambda$ are precisely those in $\Sigma$.
			Hence, as a set system, the faces of $\Phi$ are those in $\left(\Delta_{a} \cup \Delta_{b}\right) \setminus \Sigma$.
			The condition on $\Sigma$ in \labelcref{con:partition} makes $\Phi$ an actual relative simplicial complex.
		\end{proof}

		\begin{remark}
			Unless $\Delta_{a}$ or $\Delta_{b}$ is the empty complex, the intersection $\Delta_{a} \cap \Delta_{b}$ is never empty since $\varnothing$ is always a face of a nonempty simplicial complex.
			Hence, \Cref{stt:gluing_lemma} \labelcref{con:intersection_embedded} forces at least one of the complexes $\Phi_{a}$ or $\Phi_{b}$ to be a proper relative simplicial complex.
			\label{stt:intersection_deltas_non_empty}
		\end{remark}

		Since the partitioning schemes of the relative simplicial complexes of $\Phi_{a}$ and $\Phi_{b}$ do not mingle with each other, it is immediate to see the following:

		\begin{corollary}
			If $\Phi_{a}, \Phi_{b}$ and $\Phi$ as in \Cref{stt:gluing_lemma} are pure and of the same dimension, then
			$h(\Phi) = h(\Phi_{a}) + h(\Phi_{b})$.
			\label{stt:adding_up_h-vectors}
		\end{corollary}

		We will often use a particular case of \Cref{stt:gluing_lemma}.
		For the shelling portion of the statement, see the previous work of the author with Woodroofe in \cite[Lemma 2.3]{Santamaria/Woodroofe:2021}.
		
		\begin{lemma}[Shelling-like Gluing Lemma]
			Let
			$\Phi_{a} = \left( \Delta_{a}, \Gamma_{a} \right)$ and
			$\Phi_{b} = \left( \Delta_{b}, \left( \Delta_{a} \cap \Delta_{b} \right) \cup \Gamma_{b} \right)$
			be partitionable complexes such that no facet of one is properly contained in a facet of the other, where $\Gamma_{b}$ is a subcomplex of $\Delta_{b}$.
			If
			$\Delta_{a} \cap \Gamma_{b} \subseteq \Gamma_{a}$, the complex
			$\Phi = \left( \Delta_{a} \cup \Delta_{b}, \Gamma_{a} \cup \Gamma_{b} \right)$
			is partitionable.

			Furthermore, if $\Phi_{a}$ and $\Phi_{b}$ are shellable, then $\Phi$ is also shellable.
			\label{stt:assymetrical_gluing_lemma}
		\end{lemma}

		\begin{proof}
			The conditions in \Cref{stt:gluing_lemma} are met by $\Phi_{a}$ and $\Phi_{b}$ as follows:
			the relative part of $\Phi_{b}$ includes $\Delta_{a} \cap \Delta_{b}$ by definition, so the condition \labelcref{con:intersection_embedded} holds.
			It remains to prove that the complex $\Sigma$ in \labelcref{con:partition} (i.e.\ the relative part of $\Phi$) is actually $\Gamma_{a} \cup \Gamma_{b}$.
			To do this, we will treat the complexes as set systems as in \Cref{stt:gluing_lemma}.

			Let $\Gamma'_{b}$ be the relative part of $\Phi_{b}$, so $\Gamma'_{b} \colonequals \left( \Delta_{a} \cap \Delta_{b} \right) \cup \Gamma_{b}$, then
			\begin{equation}
				\Sigma
					\, = \,
				\left( \Gamma_{a} \setminus \Delta_{b} \right)
					\, \cup \,
				\left( \Gamma'_{b} \setminus \Delta_{a} \right)
					\, \cup \,
				\left( \Gamma_{a} \cap \Gamma'_{b} \right).
				\label{equ:union_shelling_like_lemma}
			\end{equation}
			After some set-based computations, we have that
			$\Gamma_{a} \cap \Gamma'_{b}
				\, = \,
			\left( \Delta_{b} \cap \Gamma_{a}\right)
				\, \cup \,
			\left( \Delta_{a} \cap \Gamma_{b} \right)$.
			Then, since
			$\Delta_{a} \cap \Gamma_{b} \subseteq
			 \Delta_{b} \cap \Gamma_{a} \subseteq
			 \Gamma_{a}$,
			it follows that
			$\Gamma_{a} \cap \Gamma'_{b}
				\, =  \,
			\Delta_{b} \cap \Gamma_{a}$.
			Thus, the union of the first and the last sets on the right-hand side of \labelcref{equ:union_shelling_like_lemma} yields $\Gamma_{a}$.

			On the other hand, it is easy to see that the set $\Gamma'_{b} \setminus \Delta_{a}$, appearing as the middle set in the union  \labelcref{equ:union_shelling_like_lemma}, is the same as $\Gamma_{b} \setminus \Delta_{a}$.
			Moreover, since
			$\Gamma_{b} \cap \Delta_{a} \subseteq
			 \Gamma_{a} \subseteq
			 \Sigma$,
			it holds that
			$\left( \Gamma_{b} \setminus \Delta_{a} \right)
				\, \cup \,
			 \left( \Gamma_{b} \cap \Delta_{a} \right)
			 	\, = \,
			 \Gamma_{b}
			$.
			Henceforth, $\Gamma_{a} \cup \Gamma_{b} = \Sigma$, as desired.

			Finally, the case when $\Phi_{a}$ and $\Phi_{b}$ are shellable, was proved in \cite[Lemma 2.3]{Santamaria/Woodroofe:2021} (see also \Cref{stt:erratum}).
		\end{proof}

		\begin{remark}
			To ensure shellability, \cite[Lemma 2.3]{Santamaria/Woodroofe:2021} requires the complex $\Gamma'_{b} = \left( \Delta_{a} \cap \Delta_{b} \right) \cup \Gamma_{b}$ to be nonempty.
			The condition is redundant, since $\Gamma'_{b}$ is never empty (see \Cref{stt:intersection_deltas_non_empty}).

			Another condition that was silently assumed by \cite[Lemma 2.3]{Santamaria/Woodroofe:2021} was that the facets in the complexes $\Phi_{a}$ and $\Phi_{b}$ must be facets of $\Phi$.
			The condition appears now explicitly stated in \Cref{stt:gluing_lemma,stt:assymetrical_gluing_lemma}.
			\label{stt:erratum}
		\end{remark}

		\begin{remark}
			The partitionings constructed by \Cref{stt:assymetrical_gluing_lemma} can be seen as generalizing those yielded by shellings, as we now explain.

			Suppose that $\Delta$ has a shelling order $\sigma_{1}, \dotsc, \sigma_{m}$, and proceed inductively on $m$.
			The base case with $m=1$ is trivial.
			Now suppose
			$\Delta_{t-1} \colonequals \langle \sigma_{1}, \dotsc, \sigma_{t-1} \rangle$ partitionable, and take
			$\Phi_{a} = \left(\Delta_{t-1}, \varnothing \right)$
			and
			$\Phi_{b} = \left(\langle \sigma_{t} \rangle, \langle \sigma_{t} \rangle \cap \Delta_{t-1} \right)$.
			Notice that the elements of $\Phi_{b}$, are precisely the same as $[R(\sigma_t), \sigma_t]$.
			Since the conditions of \Cref{stt:assymetrical_gluing_lemma} are met, $\Delta_t$ is partitionable.

			More generally, \Cref{stt:assymetrical_gluing_lemma} may be useful for finding partitioning schemes of nice enough constructible complexes.
			\label{stt:alt_shellable_are_partitionable}
		\end{remark}
	
		We believe that \Cref{stt:gluing_lemma,stt:assymetrical_gluing_lemma} may have broader use.
		Our main use in this paper will be to prove partitionability of certain $2$-dimensional spaces.
		The strategy is as follows: we split each space into more digestible ``chunks'' where partitionability is easier to analyze, then we glue them back with the help of these lemmas.
		In order to illustrate this strategy, we revisit the fact that the non-shellable Rudin's $3$-ball is partitionable \cite[\S III.2, Proposition 2.8]{Stanley:1996}.

		\begin{example}[Rudin's non-shellable $3$-ball]
			In 1958, Rudin constructed a non-shellable triangulation of the $3$-ball with a convex geometric realization \cite{Rudin:1958}.
			It is also known to be non-evasive \cite{Benedetti/Lutz:2013}.
			Rudin's ball has $f$-vector $(1,14,66,94,41)$ and $h$-vector $(1,10,30,0,0)$.
			We will prove that it is partitionable by breaking it up into smaller partitionable relative complexes
			$\left(\Delta_{a}, \varnothing \right)$
			and
			$\left(\Delta_{b}, \Delta_{a} \cap \Delta_{b} \right)$.

			We use the list of facets given in \cite{Hachimori:web}.
			Take $\{ \mathtt{1, \dotsc, 9, A, \dotsc, E} \}$ as ground set and let $\Delta_{a}$ be the shellable ball formed by the following facets (listed in a shelling order):
			\begin{center}
				\texttt{
				\begin{tabular}{ l l l l l l l l l }
					BCDE, &ABCE, &9BDE, &8CDE, &7BCD, &58CD, &3ABE, &36AE, &36AB, \\
					6ABC, &26AC, &2ACE, &24AE, &28CE, &248E, &268C, &48DE, &4ADE, \\
					458D, &458C, &68BC, &6ADE, &69DE, &48BC, &47BC, &17BD, &19BD.
				\end{tabular}}
			\end{center}
			The remaining facets in Rudin's ball generate the complex $\Delta_{b}$:
			\begin{center}
				\texttt{
				\begin{tabular}{ l l l l l l l l l}
					137D, &139D, &39CD, &59CD, &37CD, &347C, &569D, &347B, &37BE,\\
					157B, &159B, &59BE, &569E, &57BE. & & & &
				\end{tabular}}
			\end{center}
			The intersection of the two complexes is shellable, with shelling order listed below:
			\begin{align*}
				\Delta_{a} \cap \Delta_{b}
								=	\left\langle
									\mathtt{17B, 17D, 19B, 9BE, 7CD, 47B, 47C, 19D, 69D, 5CD, 69E, 3BE}
									\right\rangle.
			\end{align*}
			Notice that the facets of $\Delta_{a} \cap \Delta_{b}$ belong to the boundary of $\Delta_{a}$ and $\Delta_{b}$, and that $\Delta_{a} \cap \Delta_{b}$ is homeomorphic to a $2$-disk.

			The complex $\left(\Delta_{b}, \Delta_{a} \cap \Delta_{b} \right)$ has $h$-vector $(0,0,14,0,0)$.
			This means that if the complex is partitionable, then any partitioning scheme for $\left(\Delta_{b}, \Delta_{a} \cap \Delta_{b} \right)$ matches the fourteen edges with the fourteen facets.
			Such a partitioning is as follows:
			\begin{align*}
				\left(\Delta_{b}, \Delta_{a} \cap \Delta_{b} \right) =
				&
				\left[ \mathtt{13, 137D} \right]   \, \sqcup \,
				\left[ \mathtt{15, 159B} \right]   \, \sqcup \,
				\left[ \mathtt{34, 347B} \right]   \, \sqcup \,
				\left[ \mathtt{37, 37BE} \right]   \, \sqcup \,
				\left[ \mathtt{39, 139D} \right]   \, \sqcup \,\\
				&
				\left[ \mathtt{3C, 347C} \right]   \, \sqcup \,
				\left[ \mathtt{3D, 37CD} \right]   \, \sqcup \,
				\left[ \mathtt{56, 569D} \right]   \, \sqcup \,
				\left[ \mathtt{57, 157B} \right]   \, \sqcup \,
				\left[ \mathtt{59, 59CD} \right]   \, \sqcup \,\\
				&
				\left[ \mathtt{5B, 59BE} \right]   \, \sqcup \,
				\left[ \mathtt{5E, 569E} \right]   \, \sqcup \,
				\left[ \mathtt{7E, 57BE} \right]   \, \sqcup \,
				\left[ \mathtt{9C, 39CD} \right].
			\end{align*}
			We used computer assistance to find this partitioning.
			The problem here is computationally tractable since it is reduced to that of finding a matching in a small bipartite graph.

			The conditions of \Cref{stt:assymetrical_gluing_lemma} are met for
			$\left(\Delta_{a}, \varnothing \right)$
			and
			$\left(\Delta_{b}, \Delta_{a} \cap \Delta_{b} \right)$ with $\Gamma_{a} = \Gamma_{b} = \varnothing$, 
			hence the complex $\Delta_{a} \cup \Delta_{b}$ is partitionable, as desired.
			Finally, we notice that $\left(\Delta_{b}, \Delta_{a} \cap \Delta_{b} \right)$ is non-shellable, as otherwise the gluing would produce a shelling of Rudin's ball.
			\label{stt:rudin_ball_ex_ex}
		\end{example}

		\Cref{stt:gluing_lemma,stt:assymetrical_gluing_lemma} build up one partitionable complex out of two.
		In other situations we will instead want to modify a single partitionable complex in a manner that preserves its partitioning scheme.
		The following lemma is completely immediate.

		\begin{lemma}
			Let $\Phi$ and $\Phi'$ be relative simplicial complexes with isomorphic face posets.
			The complex $\Phi$ is partitionable if and only if $\Phi'$ is partitionable.
			\label{stt:isoposets_partitionable}
		\end{lemma}

		\begin{example}
			The following are direct consequences of \Cref{stt:isoposets_partitionable}:

			\begin{enumerate}
				\item
					\label{stt:part_minimal_representation}
					A first example of two complexes with identical posets appears when the complexes $\Phi$ and $\Phi'$ have the same minimal representation.
					As a concrete example, the leftmost and middle complexes in \Cref{fig:minimal_representations} have the same minimal representation, namely the rightmost complex.

					Observe that partitioning schemes (and also shellings) are preserved.
					This is the reason why we may consider the minimal representation of a relative complex when we are deciding partitionability (shellability) (see \Cref{stt:algorithm_feeding_rel_part}).

					\begin{figure}[!ht]
						\begin{center}
							\centering
							\includegraphics[scale=0.40]{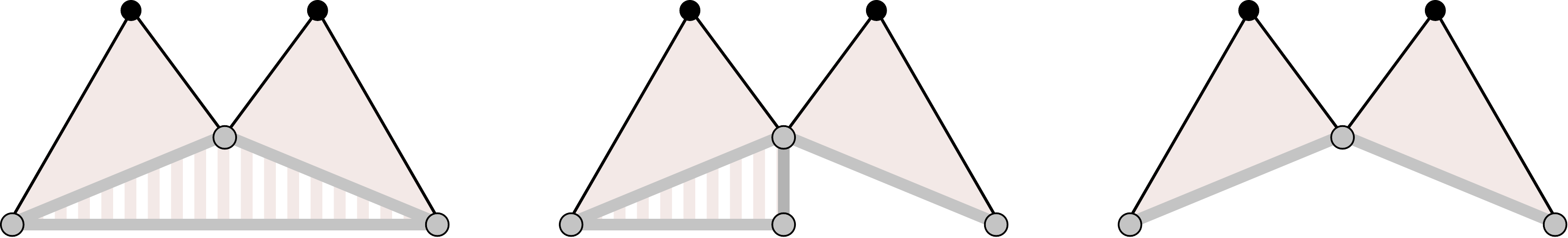}
							\caption{
								Three relative complexes with the same minimal representation.
								They share the same combinatorial closure, namely the bow-tie simplicial complex seen in \Cref{stt:bow-tie_partitionable_ex}.
								The bold light gray-colored vertices and edges, and the triangles bounded by them (left and middle), represent the relative part of each complex.}
							\label{fig:minimal_representations}
						\end{center}
					\end{figure}

				\item
					\label{stt:part_isomorphic_posets}
					More generally, \Cref{stt:isoposets_partitionable} does not depend on the dimension of the complexes involved, but strictly on their face poset structure.
					As a concrete example, \Cref{fig:isomorphic_posets} shows two partitionable simplicial complexes.
					Although the complexes do not have the same dimension, they do have isomorphic face posets.

					\begin{figure}[!ht]
						\begin{center}
							\centering
							\begin{tabular}{	>{\centering\arraybackslash}m{4.5cm}
													>{\centering\arraybackslash}m{7.8cm}
													>{\centering\arraybackslash}m{2.0cm}
													}
								\vspace{ 0.00em} \includegraphics[scale=0.40]{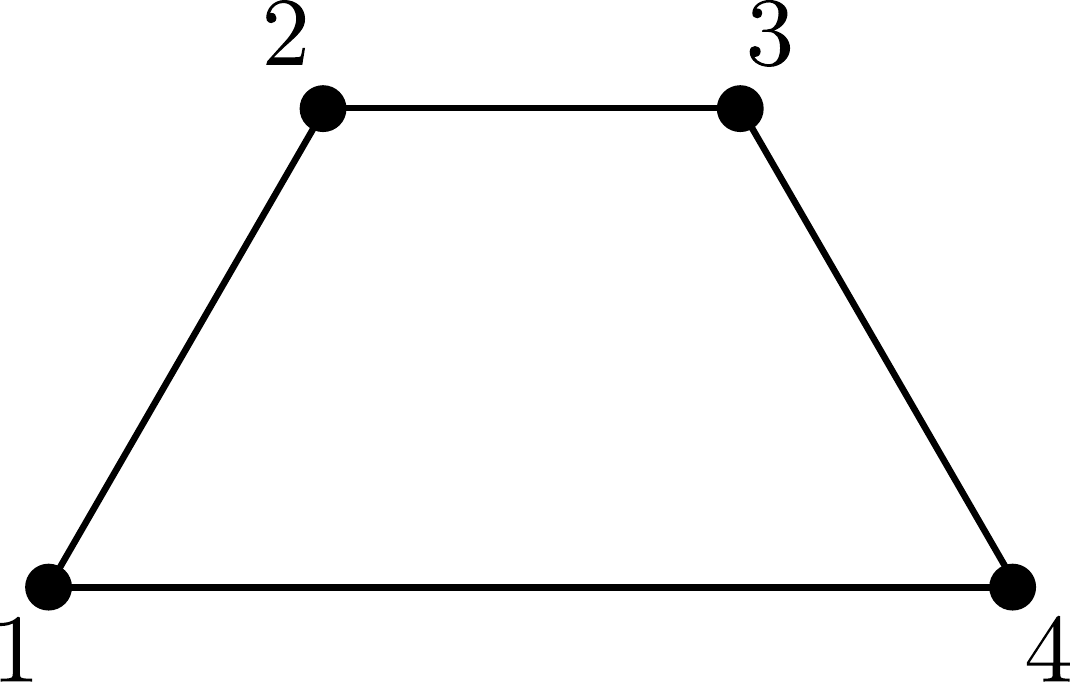} &
								\vspace{ 0.00em} \includegraphics[scale=0.46]{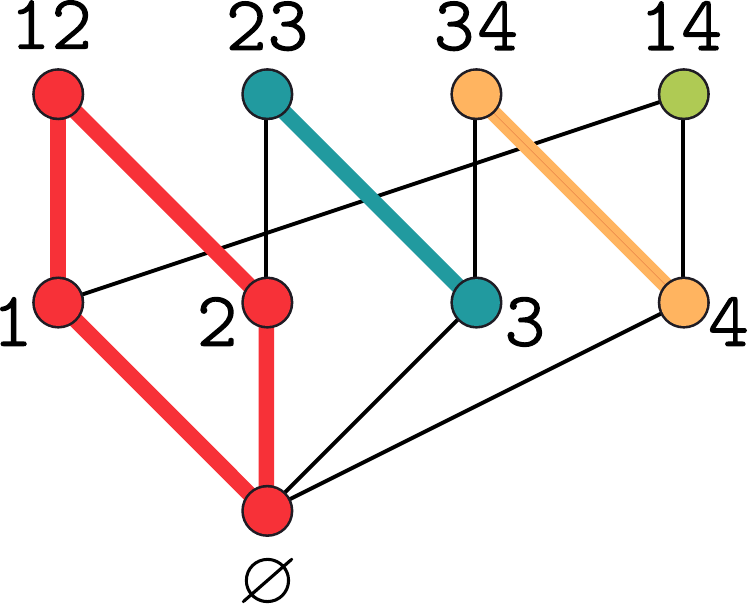} &
								\hfill
								{$\!\begin{aligned}
									\left[ \mathtt{\varnothing, 12} \right] \ & \sqcup \\
									\left[           \mathtt{3, 23} \right] \ & \sqcup \\
									\left[           \mathtt{4, 34} \right] \ & \sqcup \\
									\left[          \mathtt{14, 14} \right] \ &
								\end{aligned}$}
								\vspace{ 0.20em}
								\\
								(a) & (b) & (c) \\[1em]
								\vspace{ 0.00em} \includegraphics[scale=0.40]{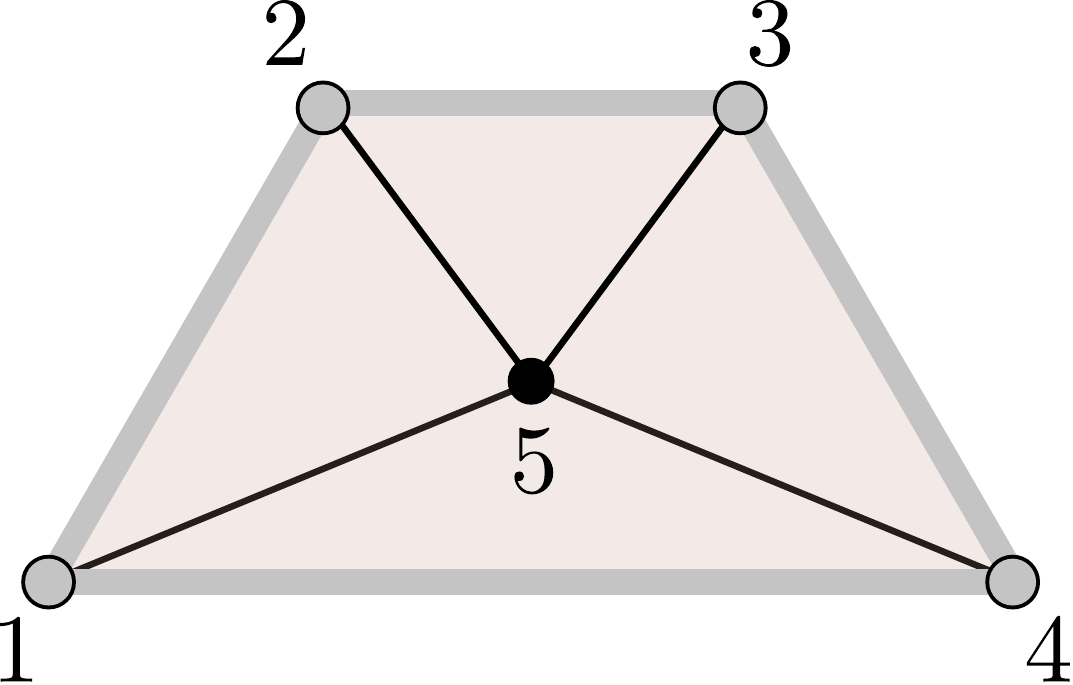} &
								\vspace{ 0.00em} \includegraphics[scale=0.46]{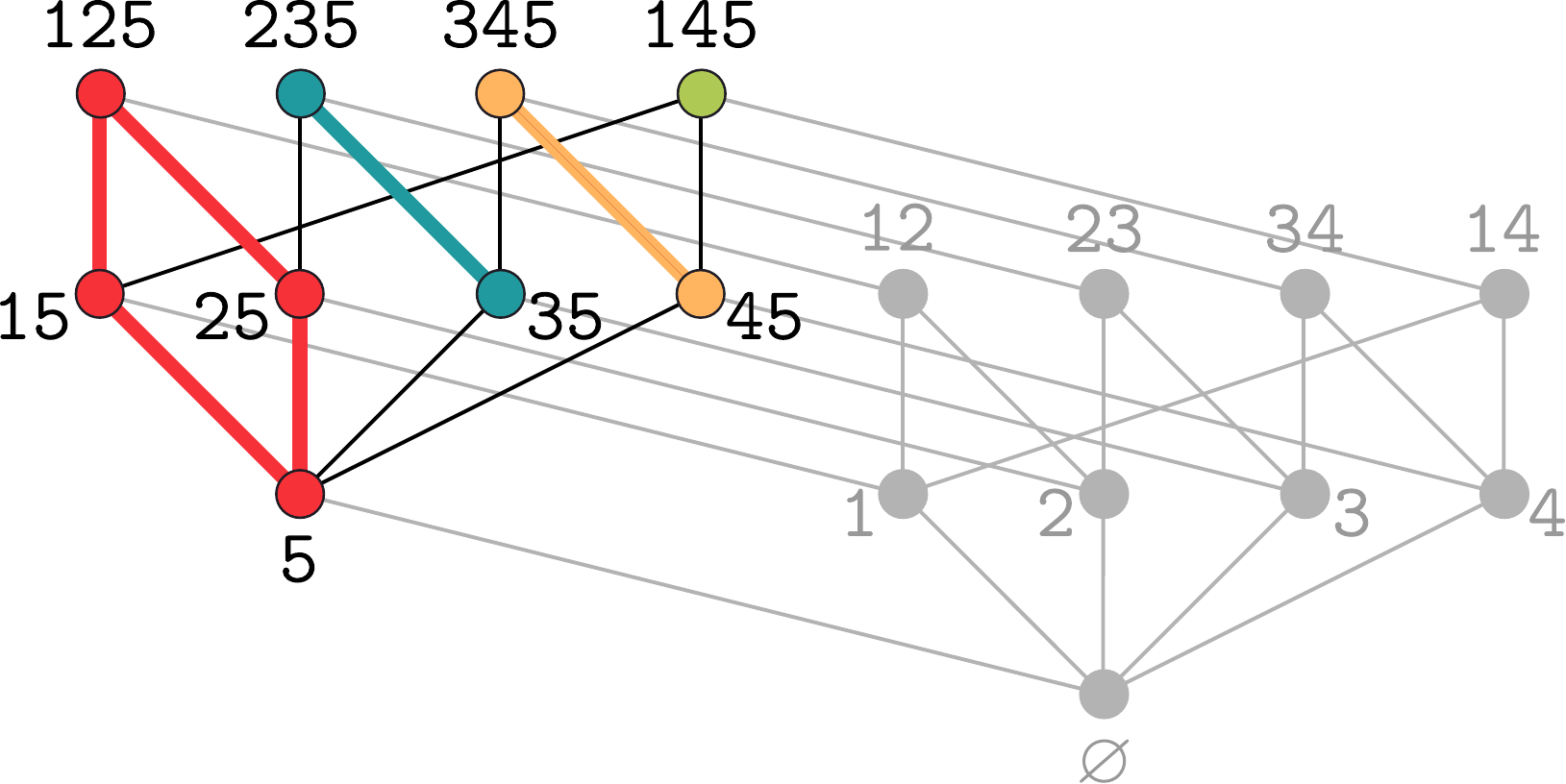} &
								\hfill
								{$\!\begin{aligned}
									\left[  \mathtt{5, 125} \right] \ & \sqcup \\
									\left[ \mathtt{35, 235} \right] \ & \sqcup \\
									\left[ \mathtt{45, 345} \right] \ & \sqcup \\
									\left[\mathtt{145, 145} \right] \ &
								\end{aligned}$}
								\vspace{ 0.20em}
								\\
								(d) & (e) & (f)
							\end{tabular}
							\caption{
								The complexes in (a) and (d) have isomorphic face posets, as shown in (b) and (e).
								The light gray-colored bold vertices and edges in the boundary of (d) form the relative part of the complex.
								Partitioning schemes are highlighted in both posets by preserving the color code of the corresponding intervals, and they are explicitly written in (c) and (f), respectively.}
							\label{fig:isomorphic_posets}
						\end{center}
					\end{figure}

			\end{enumerate}
			\label{stt:conseq_isoposets_part_ex}
		\end{example}

		We can improve \Cref{stt:isoposets_partitionable} by allowing relations to the face poset in a way that partitioning schemes are preserved.

		\begin{definition}
			A surjective dimension-preserving simplicial map $\varphi : \Delta \to \Delta'$ is a \emph{folding map} over the relative simplicial complex $\left( \Delta, \Gamma \right)$ if the following conditions hold for the map $\varphi$:

			\begin{enumerate}[label=\tt(F\arabic*)]
				\item \label{con:bijectivity}
					It induces a bijection $\left( \Delta, \Gamma \right) \leftrightarrow \left( \Delta', \Gamma' \right)$, for some $\Gamma' \subseteq \Delta'$.
				\item \label{con:facets_to_facets}
					It sends facets of $\left( \Delta, \Gamma \right)$ to facets of $\left( \Delta', \Gamma' \right)$.
			\end{enumerate}
			We say that we \emph{fold} $\left( \Delta, \Gamma \right)$ to $\left( \Delta', \Gamma' \right)$.
			\label{def:folding_map}
		\end{definition}

		We remark that the complex $\Gamma'$ in \Cref{def:folding_map} is completely determined by $\varphi$.
		Indeed, we have $ \Gamma' \colonequals \Delta' \, \setminus \, \varphi(\Delta, \Gamma)$.
		As a consequence, $\varnothing$ is a face of $\Gamma'$, except in trivial situations.
		Additionally, if $\left( \Delta, \Gamma \right)$ is pure, the conditions of \Cref{def:folding_map} lead to the equality $h(\Delta, \Gamma) = h(\Delta', \Gamma')$.

		The idea suggested by a folding map should be of a simplicial quotient map that twists and bends $\left( \Delta, \Gamma \right)$ onto itself.
		As a consequence, the image of a folding map on $\left( \Delta, \Gamma \right)$ does not necessarily preserve the topology, homotopy type, nor homology of $\Delta$ or of  $\left( \Delta, \Gamma \right)$.
		See \Cref{stt:relative_cycle_ex}.
		However, this feature is desirable for our purposes, since folding maps preserve partitioning schemes, as we now show.

		\begin{lemma}[Folding Lemma]
			Let $\left( \Delta, \Gamma \right)$ be a relative simplicial complex, and let $\varphi$ be a folding map on $\left( \Delta, \Gamma \right)$.
			If the complex $\left( \Delta, \Gamma \right)$ is partitionable, then its image under $\varphi$ is also partitionable.
			\label{stt:folding_lemma}
		\end{lemma}

		Thus, if we have a relative space where partitionability is well understood, then \Cref{stt:folding_lemma} allows us to transfer partitionings to new spaces popping up from precise manipulations of the original space.

		\begin{proof}[Proof (of \Cref{stt:folding_lemma})]
			Let $\left( \Delta', \Gamma' \right)$ be the image of $\left( \Delta, \Gamma \right)$ under $\varphi$.
			We know that $\left( \Delta', \Gamma' \right)$ is a relative simplicial complex by \labelcref{con:bijectivity}.
			Now, consider a partitioning of $\left( \Delta, \Gamma \right)$.
			We need to check that $\varphi$ carries over that partitioning from $\left( \Delta, \Gamma \right)$ to $\left( \Delta', \Gamma' \right)$.
			First, notice that $\varphi$ maps intervals into intervals.
			For every facet $\sigma \in \left( \Delta, \Gamma \right)$ we have that
			$\varphi([R(\sigma), \sigma]) = [R(\varphi(\sigma)), \varphi(\sigma)]$
			is an interval in $\Delta'$.
			This is provided by the simpliciality and dimension-preserving properties of $\varphi$.
			Also, notice that the face $\varphi(\sigma)$ is a facet of $\left( \Delta', \Gamma' \right)$ by \labelcref{con:facets_to_facets}.
			Finally, \labelcref{con:bijectivity} ensures that $[R(\varphi(\sigma)), \varphi(\sigma)]$ is an interval of $\left( \Delta', \Gamma' \right)$ and that every face in $\left( \Delta', \Gamma' \right)$ belongs to one and only one interval.
		\end{proof}

		As a cautionary remark, it is worthwhile saying that our folding maps are not related to the maps of topological pairs, where $\varphi (|\Gamma|) \subseteq |\Gamma'|$, viewed as a subspace under the subspace topology (as in \cite{Sato:1999}).
		By contrast, in a typical application of a folding map, the preimage of $\Gamma'$ is a proper subset of $\Gamma$.
		Even more, the preimage of any face in $\left( \Delta', \Gamma' \right)$ consists of one face in $\left( \Delta, \Gamma \right)$ and zero or more faces in $\Gamma$.

		The next example hints at the way we usually apply \Cref{stt:folding_lemma}.

		\begin{example}[Partitionability of $\left(\Triang{\Cycle}, \{\varnothing\} \right)$]
			Let $P_{n}$ be a path (in the graph-theoretic sense) with $n > 3$ vertices, and consider the complex $\left( P_{n}, \gSet{x} \right)$, where $x$ is a vertex of degree one.
			Let us label the vertices of $P_{n}$ like $v_1, v_2, \dotsc, v_n = x$ as we traverse from one vertex of degree one to the other.
			With exactly $n-1$ vertices and $n-1$ edges, we get a partitioning scheme for $\left( P_{n}, \gSet{x} \right)$ by creating the intervals $[v_{i}, v_{i} v_{i+1}]$, $i = 1, \dotsc, n-1$.
			By \Cref{stt:folding_lemma} we fold $\left( P_{n}, \gSet{x} \right)$ in such a way that $x$ is glued to any other vertex of choice but $v_{n-1}$ or $v_{n-2}$ (to preserve simpliciality).
			So, the obtained complex keeps the original partitioning scheme (up to some relabeling).
			In particular, when both vertices of degree one are glued to each other, the obtained complex is a cycle graph.
			Thus, any simplicial triangulation of the $1$-sphere $\Cycle$ is partitionable relative to the complex $\{\varnothing\}$.

			We observe that the face posets of $\left( P_{n}, \gSet{x} \right)$ and $\left(\Triang{\Cycle}, \{\varnothing\} \right)$ are not isomorphic despite the bijection between their facets. Hence, \Cref{stt:isoposets_partitionable} does not suffice, and we need the power of \Cref{stt:folding_lemma}.

			\label{stt:relative_cycle_ex}
		\end{example}

		\subsection{Cutting before folding}
			\label{sec:cutting_before_folding}

			As we will use the same technique later on in the paper, we want to comment on the strategy used in \Cref{stt:relative_cycle_ex}.

			We are seeking for partitionable preimages of the complex under consideration.
			In \Cref{stt:relative_cycle_ex}, we began with a cycle and selected a particular vertex $x$.
			We ``cut'' the complex at $x$ by replacing $x$ with two distinct copies of it.
			We place one copy of $x$ into the relative part, while keeping the other as a face in the relative complex.
			Each facet (edge) originally having $x$ as a vertex, now includes a distinct copy of it.
			The resulting relative complex is shellable, hence partitionable, and \Cref{stt:folding_lemma} gives the desired partitionability result.

			More generally, when we want to separate $k$ facets, we cut through the non-empty faces they have in common by making $k$ copies of each.
			One copy is kept in as a face in the relative complex, while the remaining $k-1$ copies are placed in the relative part.
			We seek to perform such cuts in a manner that yields a known shellable (or at least partitionable) relative complex.

			Once \Cref{stt:folding_lemma} is applied, it turns the complex back to its original state by identifying the replicated faces along the cuts.

		\subsection{A digression on relative partitionings of graphs}
			\label{sec:partitionable_graphs}

			The results we state here are easy to obtain and possibly well known.
			However, we present proofs both for completeness and as an example of our techniques.
			We want to generalize \Cref{stt:relative_cycle_ex}.
			First, a negative result:

			\begin{proposition}
				Let $\Delta$ be a pure $k$-dimensional simplicial complex, with $k \geq 1$.
				The complex $\left(\Delta, \{\varnothing\}\right)$ is not shellable.
				\label{stt:shell_graphs_resp_empty_face}
			\end{proposition}

			\begin{proof}
				Since the intersection of the complex induced by any facet of $\Delta$ and $\{\varnothing\}$ is not codimension one, there is no facet in $\Delta$ to start the shelling process. 
			\end{proof}

			Although $1$-dimensional complexes (graphs) are never shellable relative to $\{\varnothing\}$, many are partitionable.

			\begin{proposition}
				Let $\Delta$ be a pure $1$-dimensional simplicial complex.
				The complex $\left(\Delta, \{\varnothing\}\right)$ is partitionable if and only if no connected component of $\Delta$ is a tree.
				\label{stt:part_graphs_resp_empty_face}
			\end{proposition}

			\begin{proof}
				First, observe that a partitioning scheme of $\left(\Delta, \{\varnothing\}\right)$ must match each vertex to a unique incident edge, and the remaining edges, if any, to themselves.
				In any case, the edges have to outnumber the vertices.

					\noindent ($\Leftarrow$)
						Suppose that no connected component of $\Delta$ is a tree.
						To prove that $\left(\Delta, \{\varnothing\}\right)$ is partitionable, it is enough to see that there is a partitioning scheme for each connected component.
						Then, without loss of any generality, we may assume that $\Delta$ is a connected graph.

						Get a spanning tree $T$ of $\Delta$, and let $H$ be the subgraph generated by all the edges in $\Delta$ not included in $T$.
						See $T$ and $H$ as simplicial complexes.
						Denote as $V(H)$ the set of vertices of $H$ (seen as a $0$-dimensional complex).
						Now, select a vertex $v \in V(H)$, and consider the complexes
						$\left(T, \gSet{v} \right)$ and
						$\left(H, V(H) \setminus \{v\} \right)$.
						These two complexes meet the conditions of \Cref{stt:gluing_lemma} with $\Sigma = \{ \varnothing \}$, yielding the partitionable complex $\left(\Delta, \{\varnothing\}\right)$.
						We need to check that $\left(T, \gSet{v} \right)$ and $\left(H, V(H) \setminus \{v\} \right)$ are partitionable.
						The complex $\left(T, \gSet{v} \right)$ consists of the same number of vertices and edges, and we want to match each vertex with one of its incident edges to form intervals.
						To do that, root $T$ at $v$ and match each vertex to the edge it lies in following a bottom-up orientation (starting from the leaves, upwards).
						As for $\left(H, V(H) \setminus \{v\} \right)$, observe that it consists solely of edges and the vertex $v$.
						Match $v$ with one of the edges having it as an end-point, and match the remaining edges to themselves.

					\noindent ($\Rightarrow$)
						If $\Delta$ has a tree as a connected component, there is no way to match vertices to edges since the former set outnumbers the latter.
						\qedhere
			\end{proof}


	\section{Partitioning some relative surfaces}
		\label{sec:partitionable_surfaces}

		\subsection{Shellings on relative simplicial disks}
			\label{sec:shellable_disks}

			We are interested in partitioning $\Pplane$ and $\DunceH$.
			A main building block will be partitioning schemes of the (relative) disk.

			It is widely known that the disk (i.e.\ the $2$-ball) is shellable.
			Furthermore, every partial shelling (a disk itself) can be extended to a complete shelling of the entire disk (see, e.g.\ \cite[\S 3.5 p. 35]{Danaraj/Klee:1978a}, \cite[\S III.2 p. 84]{Stanley:1996} and \cite[\S 3.4 p. 107]{Bing:1964}).
			The following result is an easy consequence of these extendable shellings on the disk.

			\begin{theorem}
				Let $\TriangD$ be a triangulation of the $2$-dimensional disk $\Disk$ and let $\Upsilon$ be a pure connected $1$-dimensional subcomplex of its boundary.
				The complex $\left( \TriangD, \Upsilon \right)$ is shellable.
				\label{stt:rel_disk_is_partitionable}
			\end{theorem}

			This result is well known by experts in the field. 
			For completeness, and to avoid the side trip into extendable shellability, we give a proof of \Cref{stt:rel_disk_is_partitionable}.
			Our proof is essentially that of Bing in \cite[\S 3.4 Theorem 3]{Bing:1964}.

			\begin{proof}[Proof (of \Cref{stt:rel_disk_is_partitionable})]
				Proceed inductively on $m$, the number of facets of $\left( \TriangD, \Upsilon \right)$.
				The base case, when $m = 1$, is trivial.
				For the inductive step, suppose the statement holds for relative simplicial disks with less than $m$ facets, and consider a facet $\sigma$ intersecting $\Upsilon$ in a pure $1$-dimensional subcomplex.
				Such a face must exist.
				The remaining complex, say $\Phi$, can be seen as the minimal representation of $\left( \TriangD, \Upsilon \cup \gSet{\sigma} \right)$ (see \Cref{stt:conseq_isoposets_part_ex} (\labelcref{stt:part_minimal_representation})).
				Two cases may appear:

					$\bullet$ \emph{Case 1}.
						If $\sigma$ has two edges lying in the interior, and the vertex they share belongs to the boundary, $\Phi$ is split by $\sigma$ into two disks.
						Call them $\TriangDi{1}$ and $\TriangDi{2}$, and let $v$ be the vertex they hinge at.
						Let $uv$ be the edge of $\sigma$ shared with the former disk and $vw$ the one shared with the latter.
						Hence, the complexes
						$\Phi_{1} = \left(
											\TriangDi{1},
											\gSet{uv}
												\cup
											\left(\Upsilon \cap \TriangDi{1} \right)
										\right)$
						and
						$\Phi_{2} = \left(
											\TriangDi{2},
											\gSet{vw}
												\cup
											\left(\Upsilon \cap \TriangDi{2} \right)
										\right)$
						meet the inductive hypothesis and they are shellable.
						Now, build up the shelling of $\left( \TriangD, \Upsilon \right)$ by means of \Cref{stt:assymetrical_gluing_lemma}: first glue
						$\left( \gSet{\sigma}, \gSet{\sigma} \cap  \Upsilon \right) = \left( \gSet{\sigma}, \gSet{uw} \right)$
						with $\Phi_{1}$, and then glue $\Phi_{2}$ to it.

					$\bullet$ \emph{Case 2}.
						If $\sigma$ does not split the complex, then $\Phi$ ends up as a smaller relative simplicial disk meeting the condition at the boundary.
						Start the shelling with $\sigma$ and extend it through the inductive process.
			\end{proof}

			We comment that similar results to \Cref{stt:rel_disk_is_partitionable} may be given for relative cellular complexes in the sense of \cite[\S VI.6]{Bing:1983}.

		\subsection{Partitionable relative triangulations of the Möbius strip}
			\label{sec:partitioning_mobius}

			It is known that no triangulation of the Möbius strip $\Mobius$ is shellable, since it is a $2$-dimensional space homotopy equivalent to the $1$-sphere $\Cycle$.
			Nonetheless, $\Mobius$ is partitionable relative to certain subcomplexes.
			This will be a useful ingredient in the proof of our main result.

			\begin{figure}[!ht]
				\begin{center}
					\centering
					\begin{tabular}{	>{\centering\arraybackslash}b{3.6cm}
						>{\centering\arraybackslash}b{3.6cm}
						>{\centering\arraybackslash}b{3.6cm} |
						>{\centering\arraybackslash}b{3.6cm}
						}
						\includegraphics[scale=0.40]{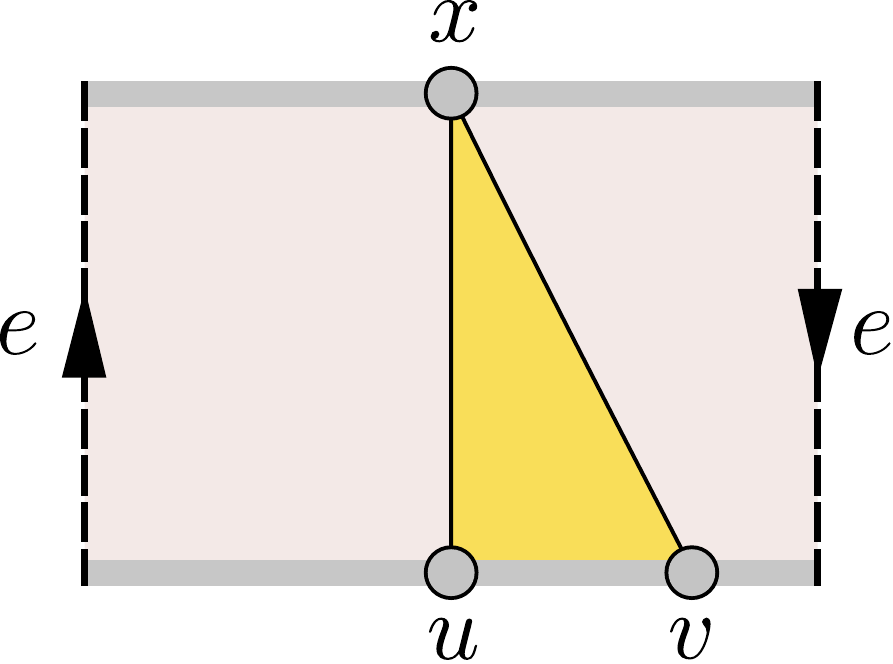} &
						\includegraphics[scale=0.40]{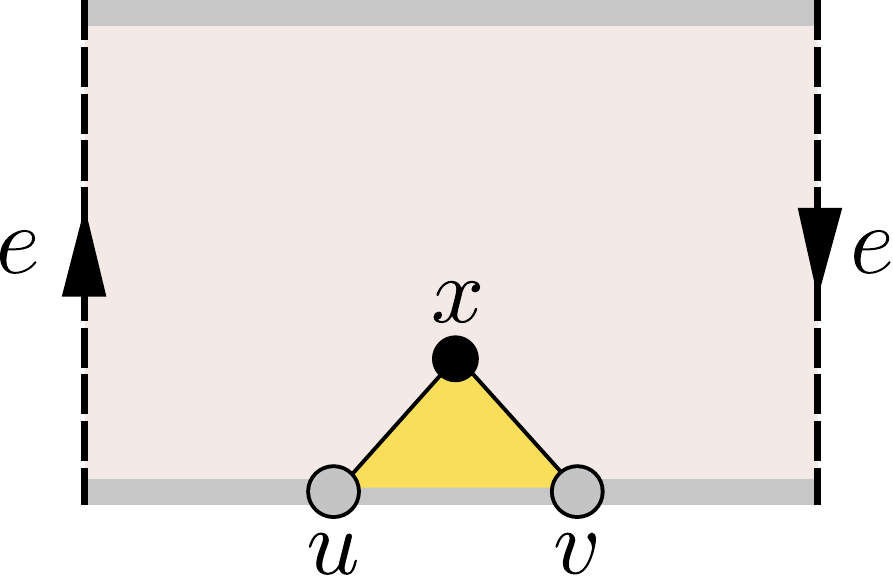} &
						\includegraphics[scale=0.40]{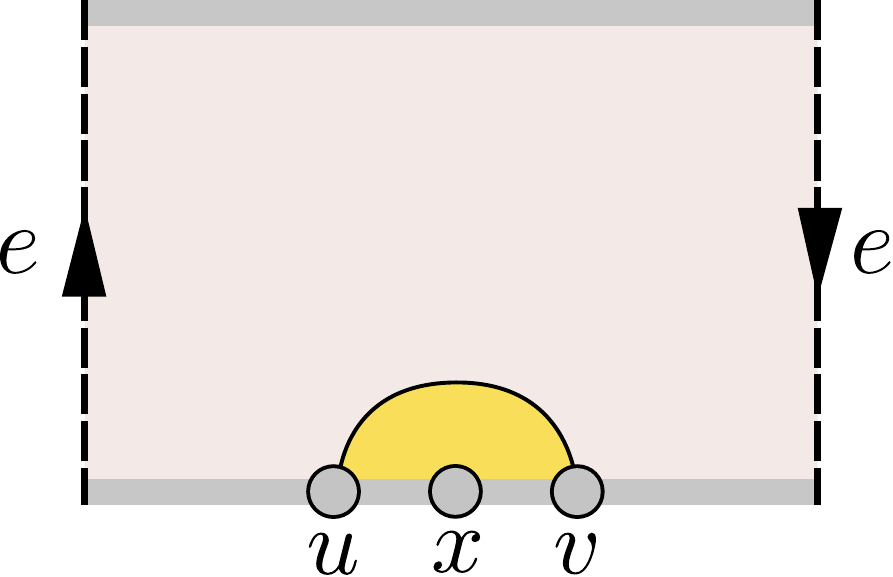} &
						\includegraphics[scale=0.40]{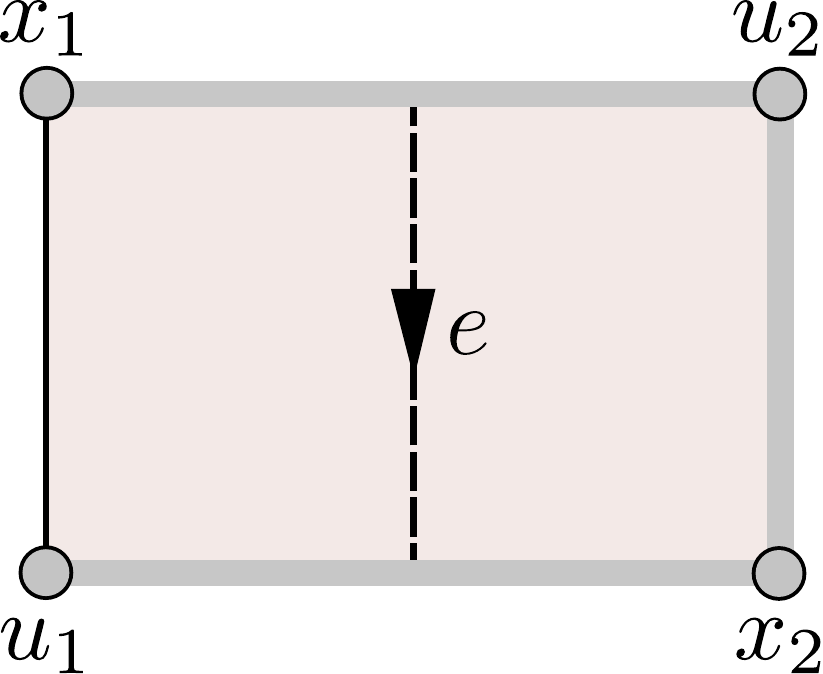} \\
						(a) & (b) & (c) & (d) \\[1em]
						\includegraphics[scale=0.40]{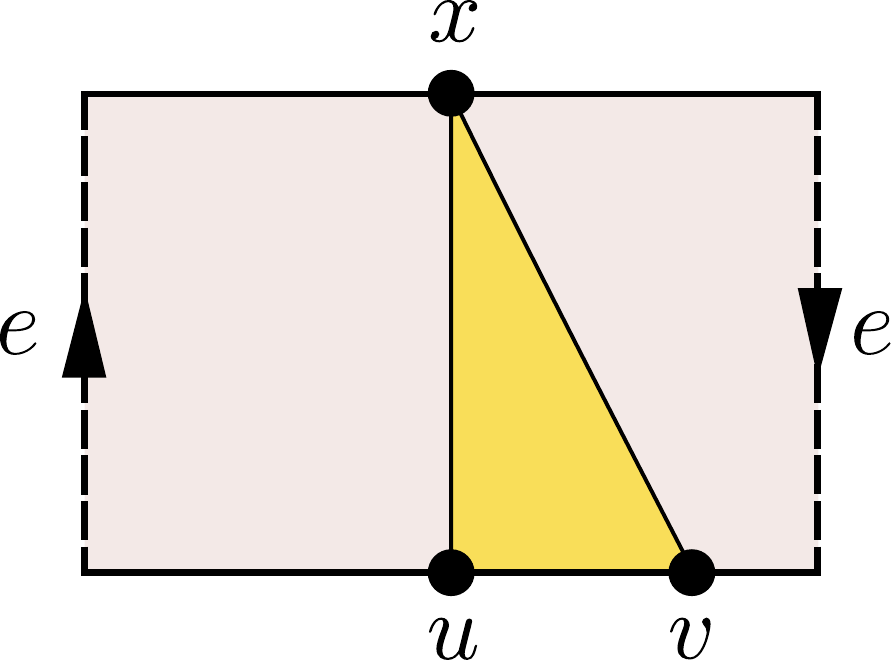} &
						\includegraphics[scale=0.40]{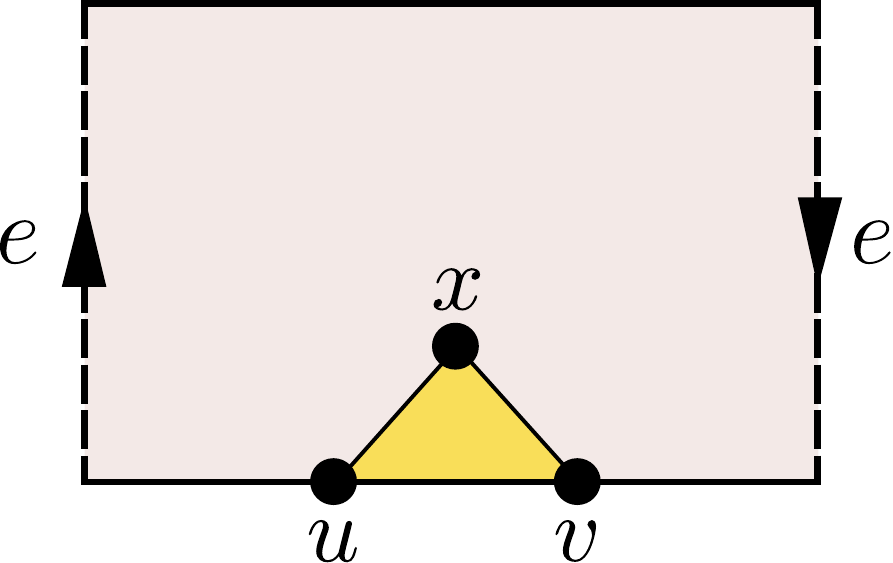} &
						\includegraphics[scale=0.40]{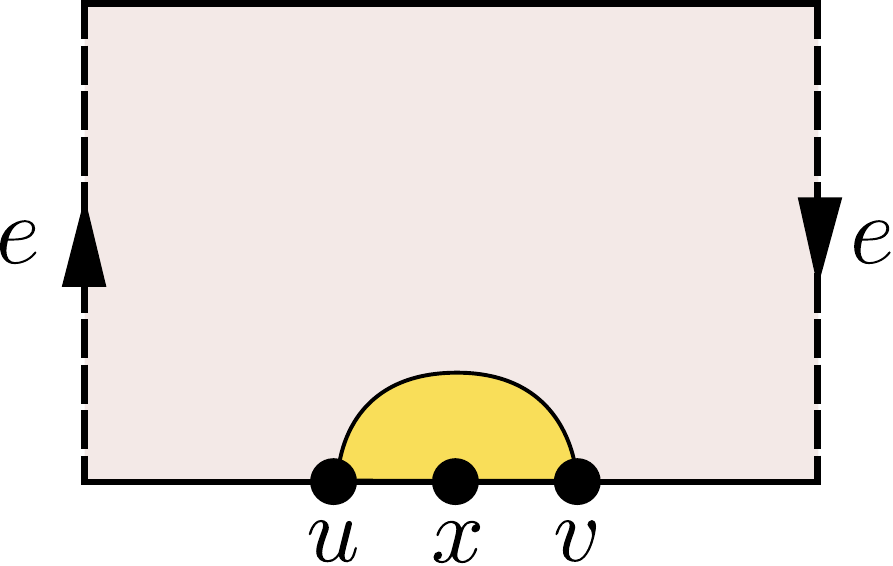} &
						\includegraphics[scale=0.40]{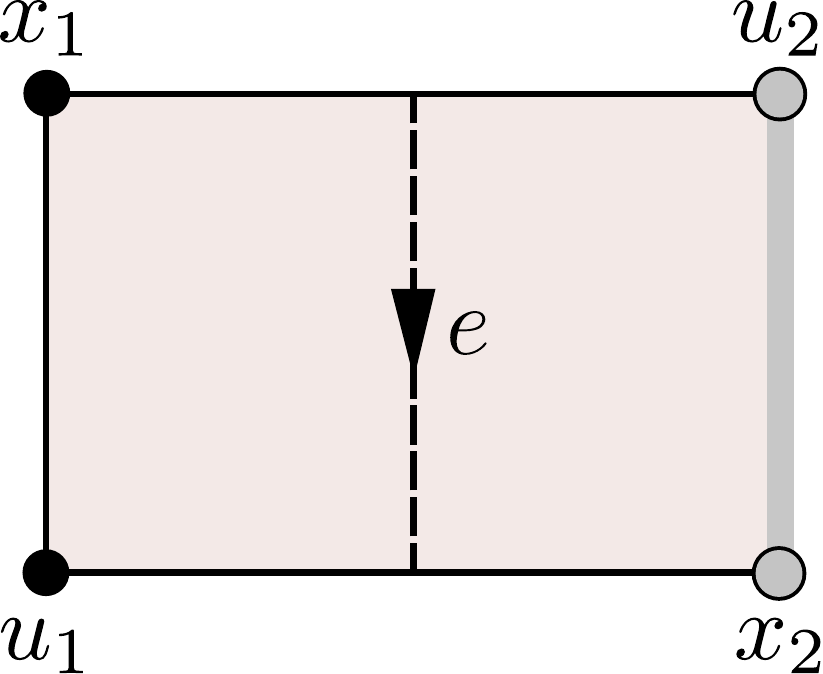} \\
						(f) & (g) & (h) & (i)
					\end{tabular}
					\caption{
						The diagrams from (a) to (c) represent triangulations of the Möbius strip $\Mobius$ relative to $\dMobius$, whereas those from (f) to (h) represent triangulations of $\Mobius$ relative to the complex $\{\varnothing\}$.
						We take out one triangle with at least one edge lying in the boundary and depict the three possible outcomes.
						In the case of (a) and (f) we cut through $ux$ to get disks relative to a portion of their respective boundaries, as it is shown in (d) and (i).
					}
					\label{fig:folding_mobius_strip}
				\end{center}
			\end{figure}

			\begin{theorem}
				Let $\TriangMo$ be a triangulation of the Möbius strip $\Mobius$.
				Then $\TriangMo$ is partitionable relative to its entire boundary.
				\label{stt:partitionable_mobius_boundary}
			\end{theorem}

			\begin{proof}
				We apply the following inductive argument.
				Pick a triangle $\sigma = uvx$ with an edge lying on the boundary of $\TriangMo$.
				Let $\dTriangMo$ be the boundary of $\TriangMo$.

				We have three cases, all of them are drawn in \Cref{fig:folding_mobius_strip} (a) to (c).
				If a case like \Cref{fig:folding_mobius_strip} (a) appears, namely, two edges in the interior and all vertices in the boundary, we can split the strip to turn it into a disk:
				take the edge $ux$ in the interior, and cut through it as explained in \Cref{sec:cutting_before_folding} to get a disk relative to the path $x_1 u_2 x_2 u_1$ (see \Cref{fig:folding_mobius_strip} (d)).
				Now, get a partitioning with \Cref{stt:rel_disk_is_partitionable} and fold back with \Cref{stt:folding_lemma}.
				Use this as the base case for the induction.

				In the case that such a triangle does not exist, consider the remaining cases (\Cref{fig:folding_mobius_strip} (b) and (c)).
				Taking the triangle $\sigma$ out of $\TriangMo$ does not split the strip and preserves its topology.
				Thus, we get a new triangulation of $\Mobius$ relative to its boundary, call it $\left(\TriangMo', \dTriangMo' \right)$.
				More precisely, such a triangulation is better described as the minimal representation of $\left( \TriangMo, \dTriangMo \cup \gSet{\sigma} \right)$.
				By the inductive hypothesis, $\left(\TriangMo', \dTriangMo' \right)$ is partitionable.
				Use \Cref{stt:assymetrical_gluing_lemma} to glue back the relative triangle we took out from the original complex, i.e.\
				$\left( \gSet{\sigma}, \gSet{\sigma} \cap \dTriangMo \right)$.
				This yields a partitioning scheme of $\left(\TriangMo, \dTriangMo \right)$.
				Observe that \Cref{stt:assymetrical_gluing_lemma} was invoked here with 
				$\Phi_{a} = \left( \gSet{\sigma}, \gSet{\sigma} \cap \dTriangMo \right)$, 
				and $\Gamma_{b} = \gSet{\dTriangMo \setminus \gSet{\sigma}}$.
			\end{proof}

			For the sake of completeness, we present an analogous result to that of \Cref{stt:partitionable_mobius_boundary}.
			Use \Cref{fig:folding_mobius_strip} (f) to (i) as a pictorial reference.

			\begin{proposition}
				Let $\TriangMo$ be a triangulation of the Möbius strip $\Mobius$.
				Then $\TriangMo$ is partitionable relative to the complex $\{\varnothing\}$.
				\label{stt:partitionable_mobius_empty_face}
			\end{proposition}

			\begin{proof}
				The proof is completely analogous to that of \Cref{stt:partitionable_mobius_boundary}.
				However, for the base case, the complex we get is a disk relative to a single edge (see \Cref{fig:folding_mobius_strip} (i)).
				Apply \Cref{stt:rel_disk_is_partitionable} and \Cref{stt:folding_lemma}.
				As for the inductive step, with $uvx$ as before, we must consider the complexes
				$\left(\TriangMo \setminus \left[ uv, uvx \right], \{ \varnothing \} \right)$
				and
				$\left( \gSet{uvx}, \gSet{ux, vx} \right)$
				for the case illustrated in \Cref{fig:folding_mobius_strip} (g),
				and the complexes
				$\left(\TriangMo \setminus \left[ x, uvx \right], \{ \varnothing \} \right)$
				and
				$\left( \gSet{uvx}, \gSet{uv} \right)$
				for the case in \Cref{fig:folding_mobius_strip} (h).
				Glue with \Cref{stt:assymetrical_gluing_lemma}.
			\end{proof}

			By using the same method, we can obtain analogous results for the annulus $\Annulus$.

			\begin{proposition}
				Let $\TriangAn$ be a triangulation of the annulus $\Annulus$.
				Then the following relative complexes are partitionable.
				\begin{enumerate}
					\item $\TriangAn$ relative to $\{\varnothing\}$.
					\item $\TriangAn$ relative to its entire boundary.
					\item $\TriangAn$ relative to one of the cycles in its boundary.
				\end{enumerate}
				\label{stt:partitionable_annulus}
			\end{proposition}

			\begin{proof}
				Start with a simplicial triangulation $\TriangAn$ of $\Annulus$.
				Then get the closest path between the two contours of $\TriangAn$ on its $1$-skeleton and cut through it.
				Now, for each of the three introduced cases, get a partitioning scheme with \Cref{stt:rel_disk_is_partitionable} and apply the folding argument of \Cref{stt:folding_lemma}, as in \Cref{stt:partitionable_mobius_boundary} and \Cref{stt:partitionable_mobius_empty_face}.
			\end{proof}

			We are sure that the reader can apply our cut-and-fold techniques to obtain results on partitionability of other surfaces, like the sphere, the torus, or the Klein bottle, in relative terms.


	\section{Proof of \Cref{stt:main_result}}
		\label{sec:main_theorem}

		\subsection{Partitionable triangulations of the projective plane}
			\label{sec:partitioning_pplane}

			The real projective plane $\Pplane$ is a non-orientable surface with Euler characteristic $0$.
			No triangulation of $\Pplane$ is shellable.
			In this section we prove that $\Pplane$ is partitionable.

			\begin{figure}[!ht]
				\begin{center}
					\centering
					\begin{tabular}{	>{\centering\arraybackslash}m{5.3cm}
											>{\centering\arraybackslash}m{5.3cm}}
						\vspace{ 0.0em} \includegraphics[scale=0.40]{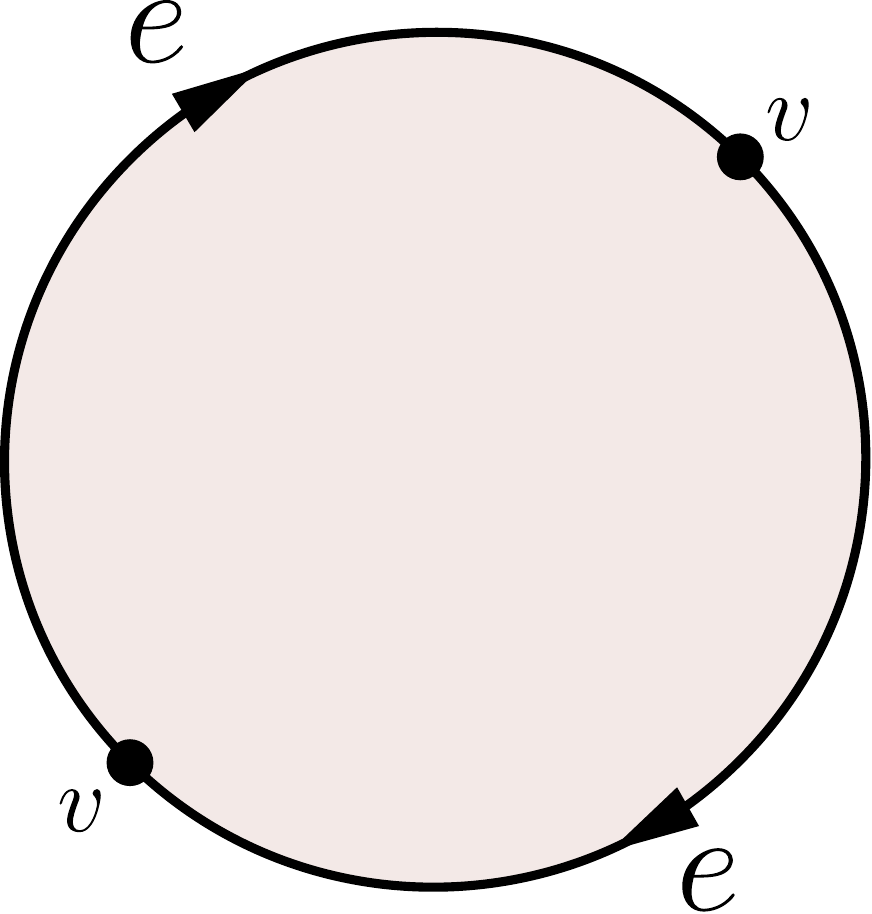} &
						\vspace{ 0.0em} \includegraphics[scale=0.40]{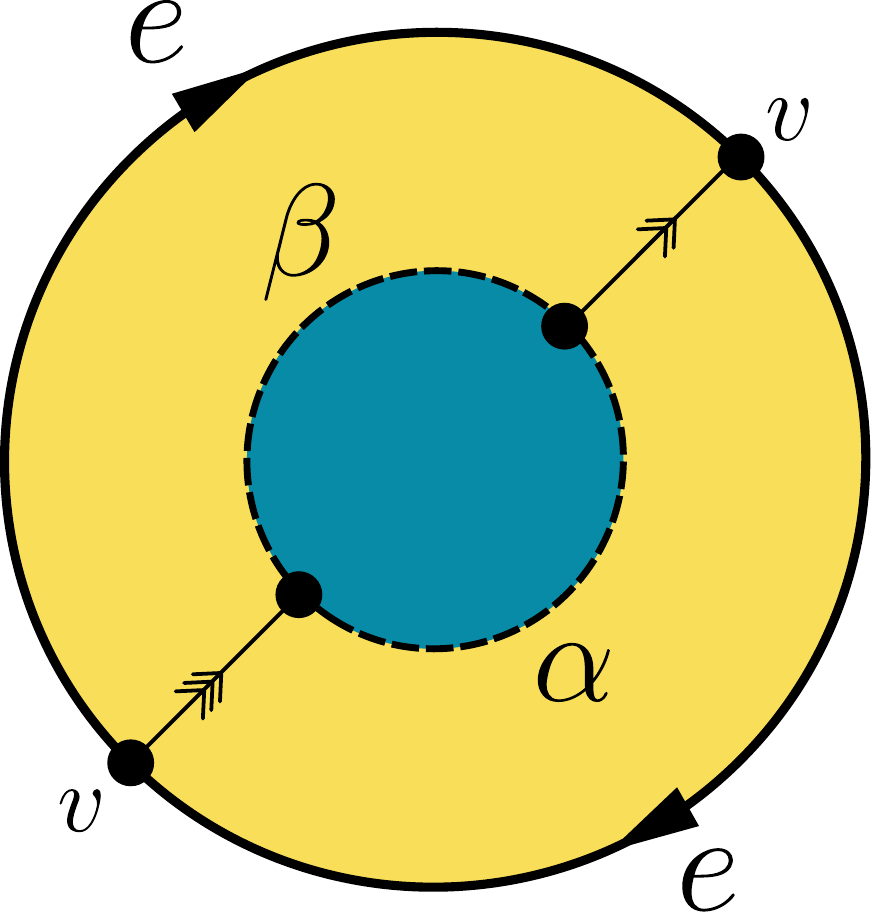} \\
						(a) & (b) \\[1em]
						\vspace{ 0.0em} \includegraphics[scale=0.40]{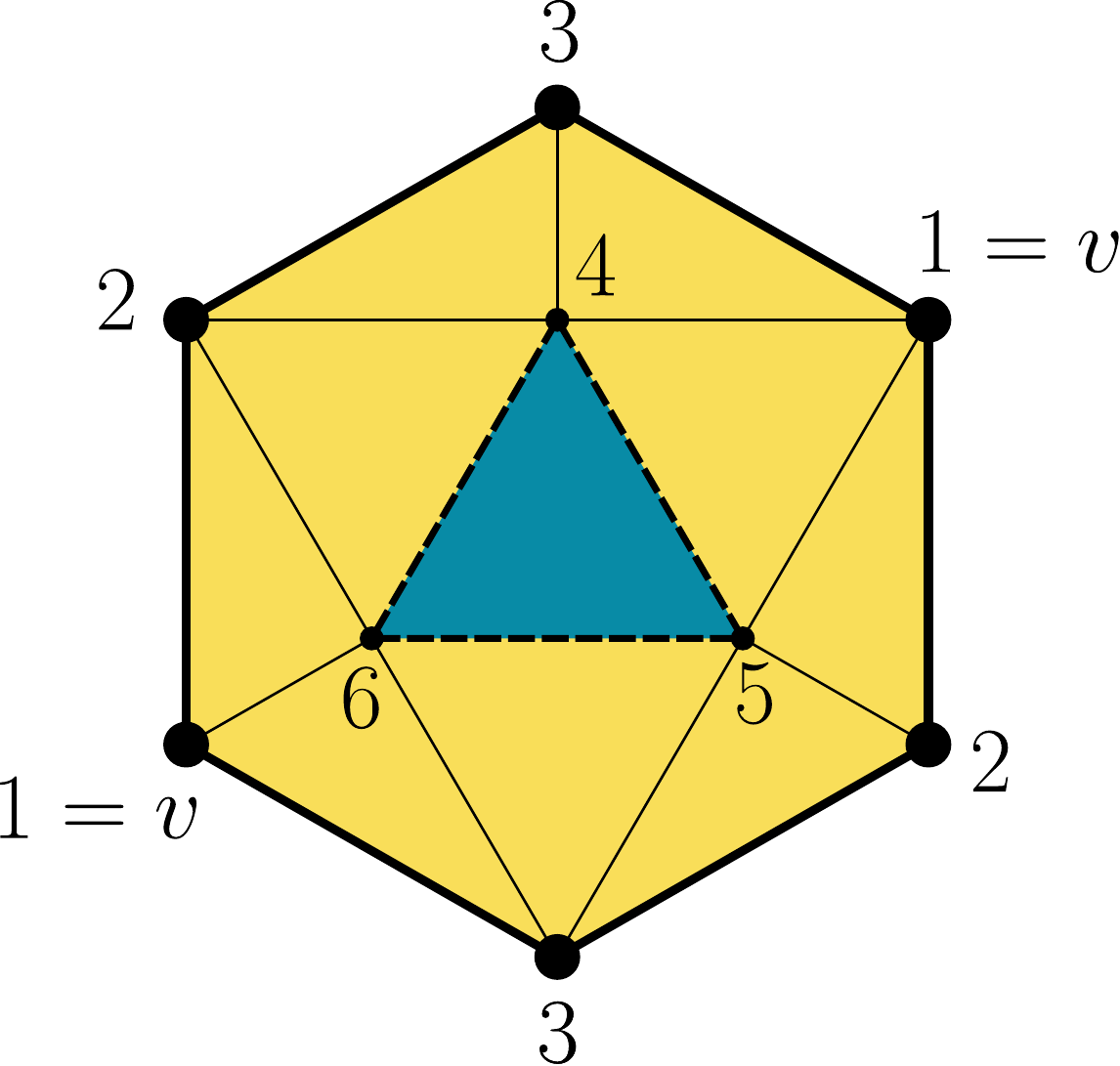} &
						\vspace{ 0.0em} \includegraphics[scale=0.40]{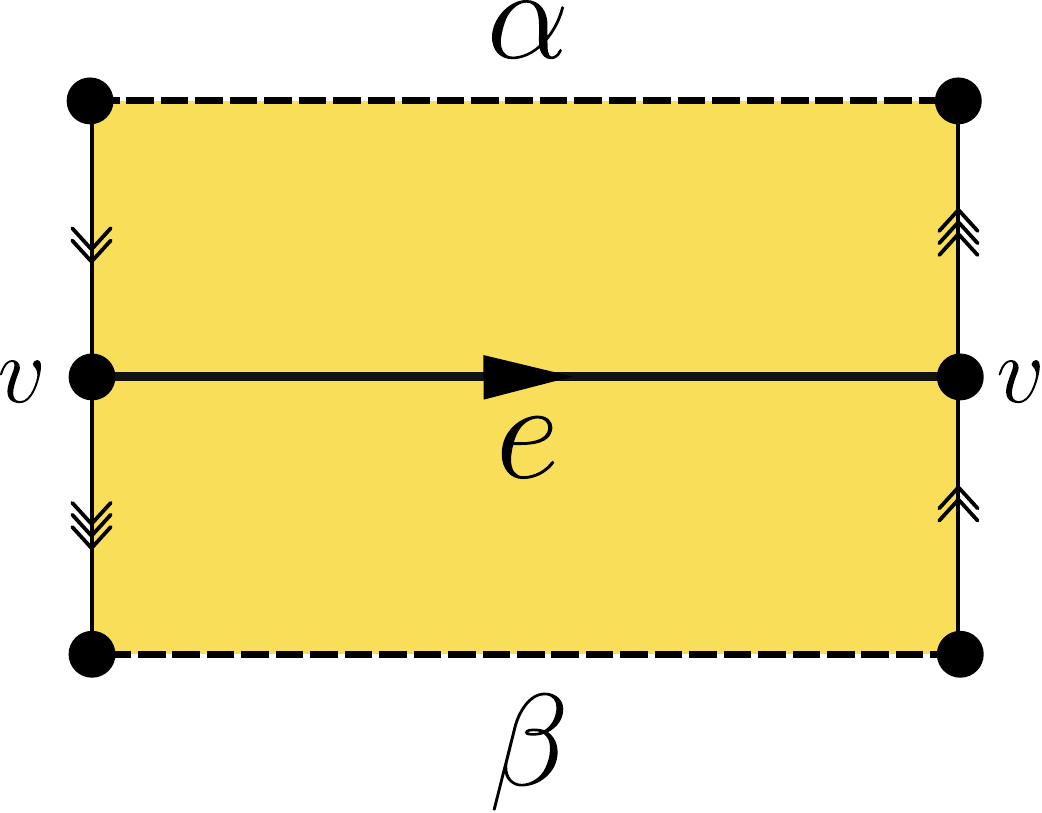} \\
						(c) & (d)
					\end{tabular}
					\caption{
						The real projective plane $\Pplane$ as a CW-complex (a).
						Sides labeled as $e$ must be identified.
						A decomposition of $\Pplane$ into a Möbius strip (yellow) and a disk (dark blue) (b). 
						This decomposition is shown in a particular triangulation of $\Pplane$ (c).
						We redraw in (d) the yellow-colored space in (b) to make clear it is a Möbius strip.}
					\label{fig:projective_plane}
				\end{center}
			\end{figure}

			\subsubsection{Decomposition of $\Pplane$}
				\label{sec:decomposition_pplane}

				We want to split $\Pplane$ into suitable partitionable subcomplexes, with the aim to glue them back by means of \Cref{stt:assymetrical_gluing_lemma}.

				\begin{lemma}
					Let $\TriangPp$ be a triangulation of the projective plane $\Pplane$, and $\sigma$ a facet of $\TriangPp$.
					Then the complex $\TriangPp \setminus \sigma$ is a triangulation of the Möbius strip $\Mobius$.
					\label{stt:decomposition_pplane}
				\end{lemma}

				\begin{proof}
					Take an arbitrary facet $\sigma$ out of $\TriangPp$.
					The complex $\gSet{\sigma}$ is obviously a simplicial disk, and its boundary is homeomorphic to a cycle $\Cycle$.
					Notice that the remaining complex $\Delta \colonequals \TriangPp \setminus \sigma$ shares its entire boundary with $\gSet{\sigma}$.
					Now, we want to recognize the space that $\Delta$ triangulates.
					Since the geometric realization of $\gSet{\sigma}$ (homeomorphic to a disk) and the space $\Pplane$ wherein it lies are surfaces, then \cite[\S 3, Corollary 3.14$_{n+1}$]{Rourke/Sanderson:1972} ensures that $|\Delta| \cong \cl(\Pplane \setminus |\gSet{\sigma}|)$ is a surface.

					On the other hand, the deletion of $\sigma$ reduces by one the Euler characteristic of $\Pplane$, so $\tilde{\chi} (|\Delta|) = -1$.
					Then, by the classification theorem of compact surfaces (see \cite[\S 6.3, Theorem 6.2]{Gallier/Xu:2013}), the complex $\Delta$ is a triangulation of either the Möbius strip $\Mobius$ or the annulus $\Annulus$.
					As the boundary of the annulus is disconnected, the only possible option is that $|\Delta| \cong \Mobius$.
				\end{proof}

				\begin{remark}
					This decomposition of $\Pplane$ into a simplicial disk and a triangulation of a Möbius strip came as no surprise.
					As a consequence of the classification of surfaces, any compact non-orientable surface is obtained from the sphere by deleting a disk and identifying the resulting boundary cycle with the boundary of a Möbius strip.
					Also, all compact surfaces with boundary are obtained by deleting several disks from a closed surface (see \cite{Gallier/Xu:2013} for further reference).
					However, there is a subtlety we had to take care of in our decomposition: deleting an arbitrary disk is not the same as deleting well-chosen one.
					We were surprised not to find a result like \Cref{stt:decomposition_pplane} explicitly stated in the literature, although we believe it is well known by the experts in the field.
				\end{remark}

			\medskip
			\phantomsection
				\label{sec:main_result_pplane}

				We are now ready to prove the first part of \Cref{stt:main_result}.

				\begin{theorem}[\Cref{stt:main_result} for $\Pplane$]
					Any triangulation of $\Pplane$ is partitionable.
					\label{stt:pplane_partitionable}
				\end{theorem}

				\begin{proof}
					\Cref{stt:decomposition_pplane} allows us to decompose a triangulation $\TriangPp$ of the projective plane into an arbitrary facet $\sigma$ and a triangulation $\TriangMo$ of the Möbius strip $\Mobius$.

					Denote the boundary of $\TriangMo$ as $\dTriangMo$.
					Following \Cref{stt:assymetrical_gluing_lemma} (and notation) we need partitionable complexes
					$\Phi_{a} = \gSet{\sigma} = \left( \gSet{\sigma}, \varnothing \right)$ and
					$\Phi_{b} = \left( \TriangMo, \dTriangMo \right)$.
					The former complex is obviously shellable, and the latter is partitionable by \Cref{stt:partitionable_mobius_boundary}.
					To check the set-theoretic conditions of \Cref{stt:assymetrical_gluing_lemma} observe that $\TriangMo \cap \gSet{\sigma} = \dTriangMo = \gSet{\sigma} \setminus \sigma$ and $\TriangMo \cup \gSet{\sigma} = \TriangPp$.
				\end{proof}

				By using the same techniques we can easily obtain a similar relative result.

				\begin{theorem}
					Any triangulation of $\Pplane$ is partitionable relative to $\{\varnothing\}$.
					\label{stt:pplane_rel_empy_face_partitionable}
				\end{theorem}

				\begin{proof}
					Entirely analogously to the previous proof, we consider the complexes
					$\Phi_{a} = \left( \TriangMo, \{\varnothing\} \right)$ and
					$\Phi_{b} = \left( \gSet{\sigma}, \gSet{\sigma} \setminus \sigma \right)$.
					The latter is obviously shellable, and the former is partitionable by \Cref{stt:partitionable_mobius_empty_face}.
				\end{proof}

				\begin{remark}
					It might be interesting to consider partitionability of $\Pspace{k}, k > 2$, with the use of our techniques.
					However, we do not pursue this problem any further in the current work.
					\label{stt:future_work_RPk}
				\end{remark}

		\subsection{Partitionable triangulations of the dunce hat}
			\label{sec:partitioning_dunce_hat}

			\begin{figure}[!ht]
				\begin{center}
					\centering
					\begin{tabular}{	>{\centering\arraybackslash}m{5.3cm}
											>{\centering\arraybackslash}m{5.3cm}}
						\vspace{ 0.20em} \includegraphics[scale=0.40]{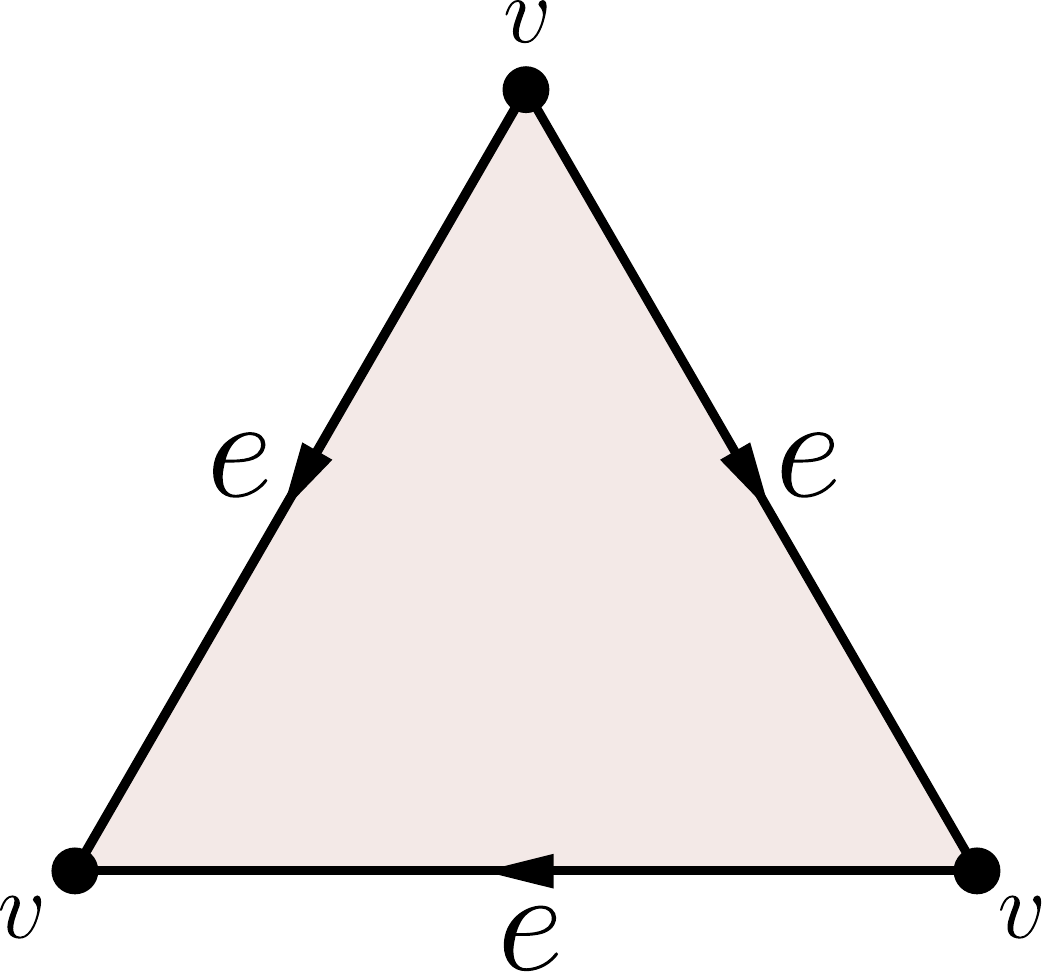} &
						\vspace{ 0.00em} \includegraphics[scale=0.40]{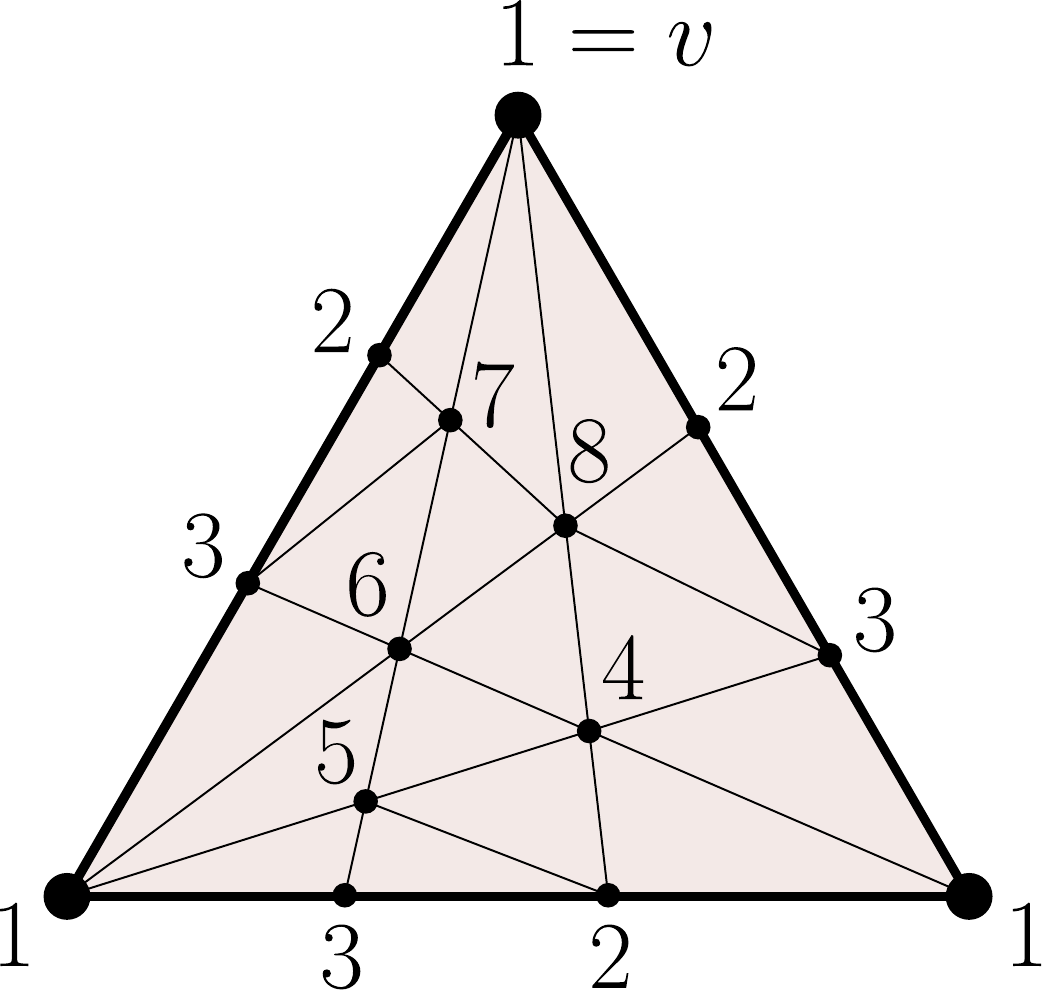} \\
						(a) & (b) \\
						\vspace{ 1.13em} \includegraphics[scale=0.40]{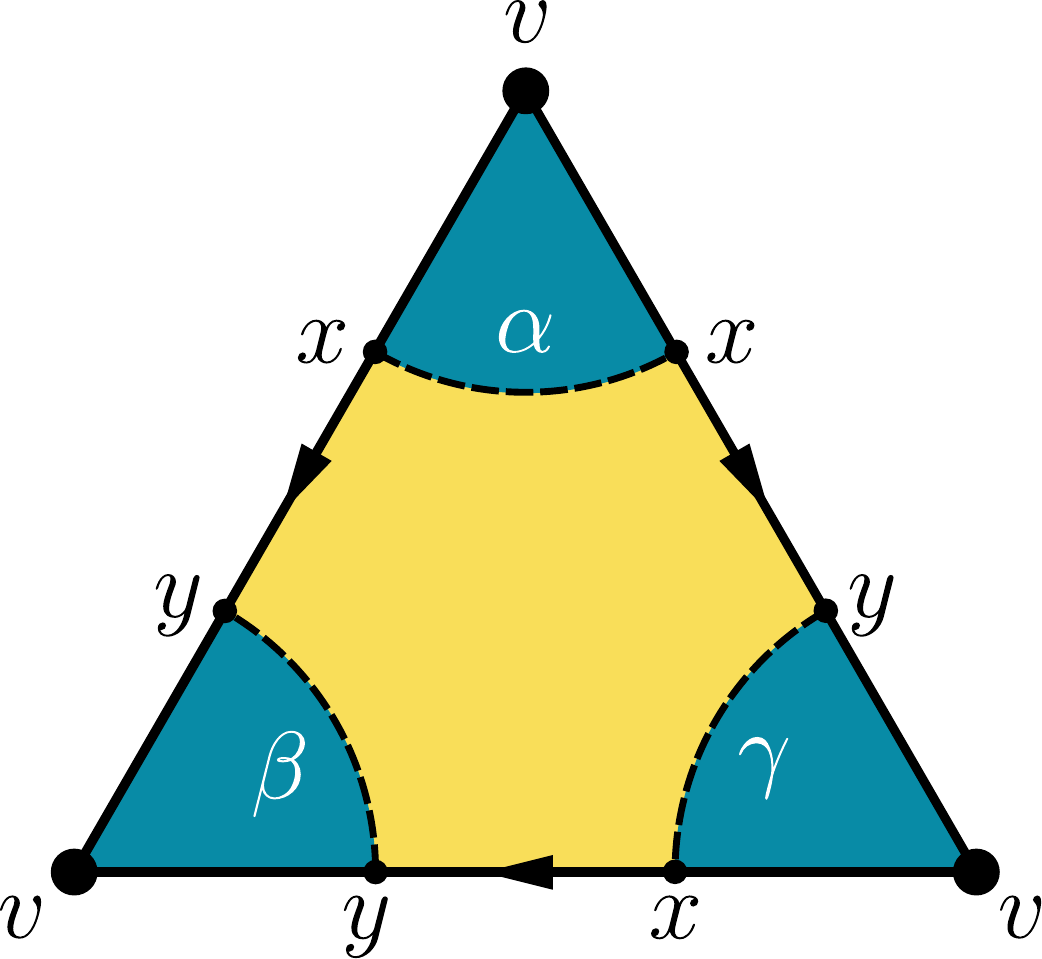} &
						\vspace{ 1.00em} \includegraphics[scale=0.40]{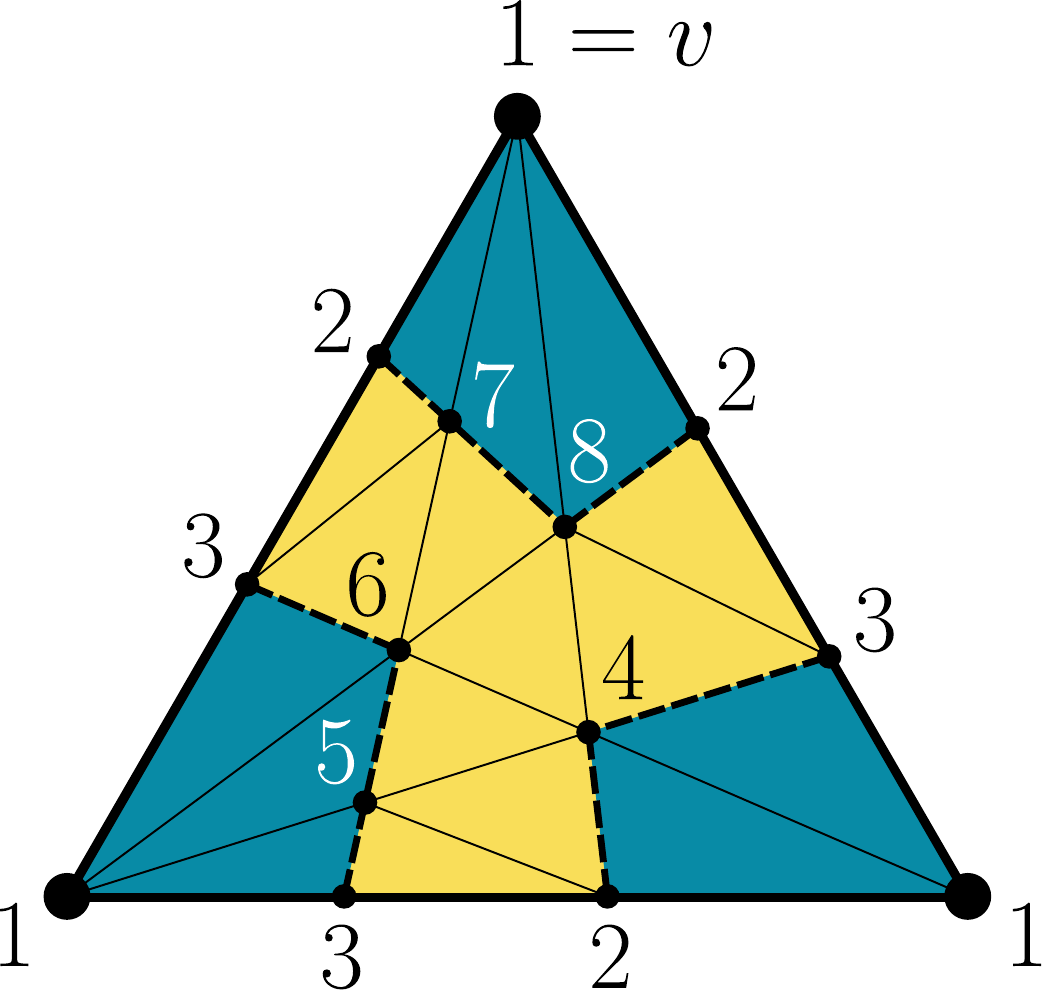} \\
						(c) & (d) \\
						\vspace{ 1.00em} \includegraphics[scale=0.40]{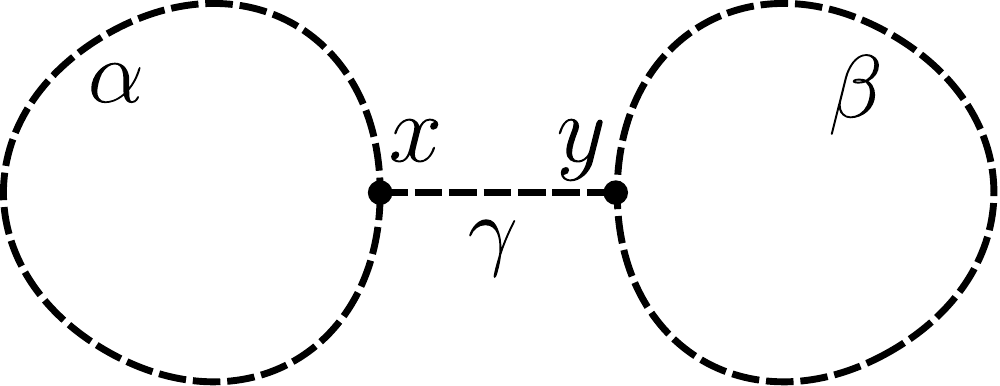} &
						\vspace{ 1.00em} \includegraphics[scale=0.40]{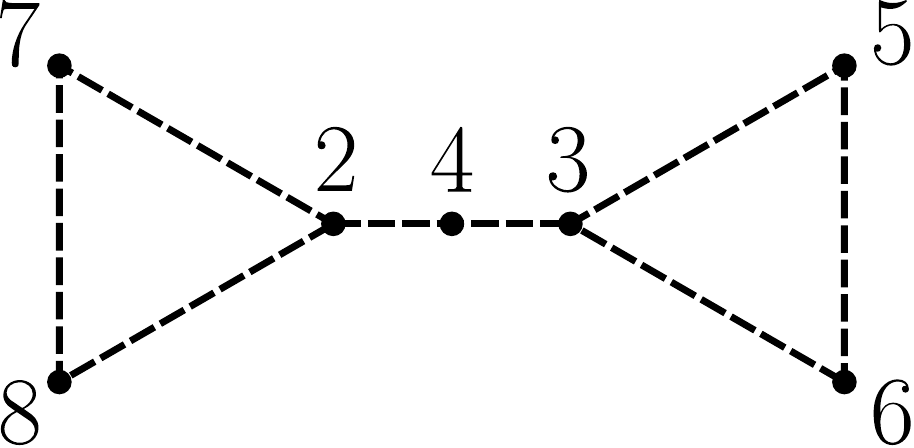} \\
						(e) & (f) \\
						\vspace{ 1.13em} \includegraphics[scale=0.40]{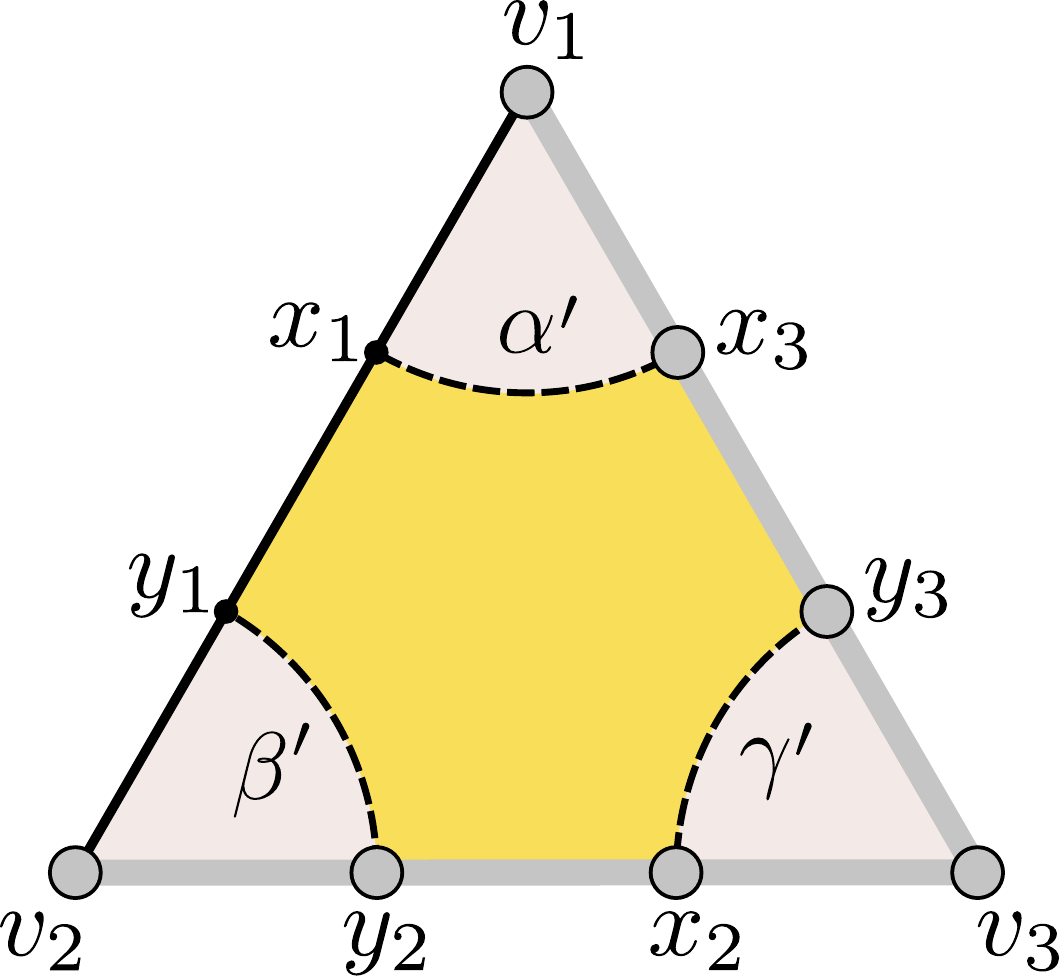} &
						\vspace{ 0.90em} \includegraphics[scale=0.40]{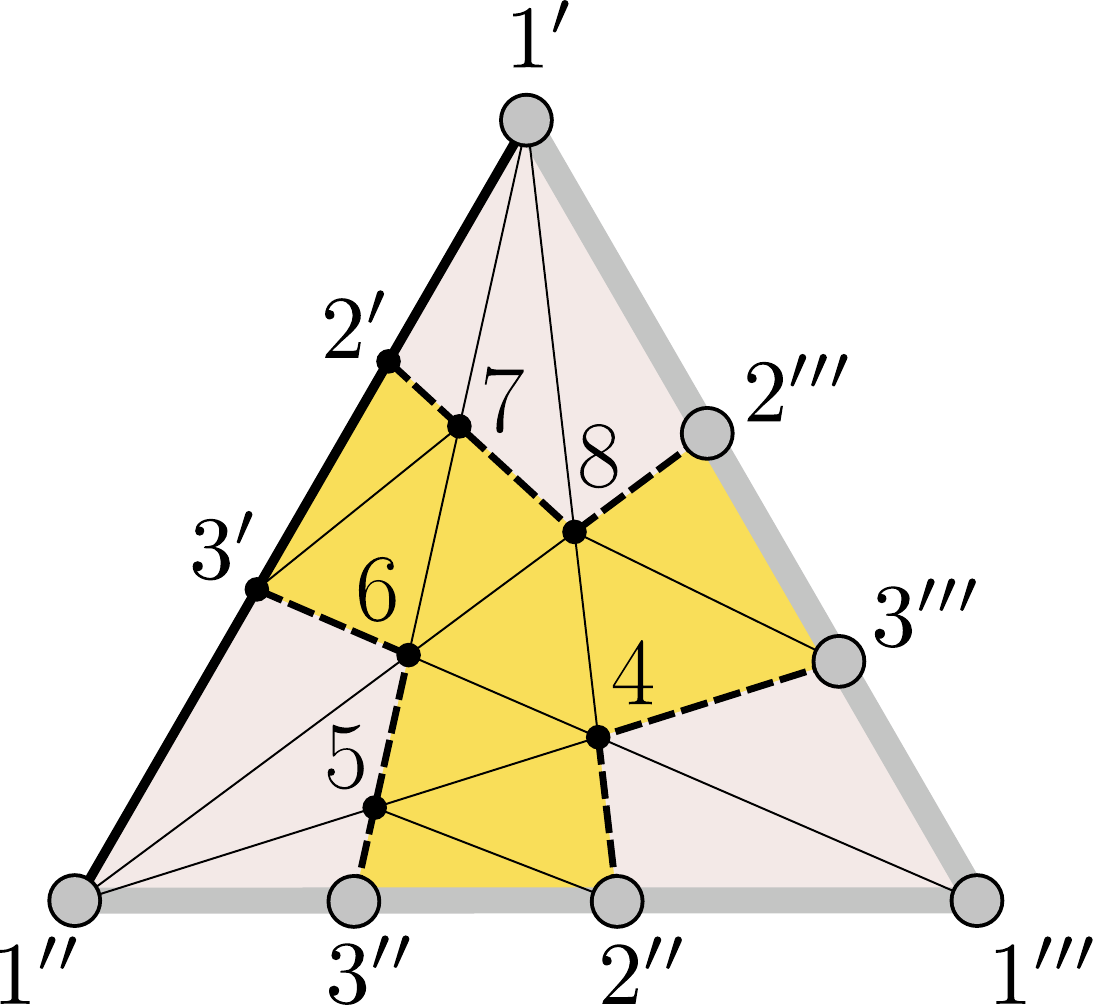} \\
						(g) & (h) 
					\end{tabular}
					\caption{
						The dunce hat $\DunceH$ as a CW-complex (a) and a triangulation $\TriangDH$ of $\DunceH$ (b).
						Sides labeled as $e$ must be identified preserving the orientation.
						In (c) and (d) the space is decomposed into the cone of $v$ (dark blue) and the deletion of $v$ (yellow).
						The boundary shared by these two spaces is shown in (e) and (f).
						To obtain a partitioning for the deletion of $v$, we first get a disk cutting along $e$, as shown in (g) and (h).
						Then we make the deletion of each copy of $v$.}
					\label{fig:dunce_hat}
				\end{center}
			\end{figure}

			The \emph{dunce hat} $\DunceH$ is the quotient space of a $2$-dimensional triangle, where the sides are identified in a non-cyclic manner.
			It was introduced by Zeeman in \cite{Zeeman:1964} (hence $\DunceH$ for the notation).
			The space $\DunceH$ is known to be contractible but non-collapsible; it is Cohen-Macaulay over any field, but it is not shellable \cite[\S III.2 p. 84]{Stanley:1996} nor even constructible \cite{Hachimori:2008}.
			In \Cref{fig:dunce_hat} we depict the space $\DunceH$ (a) and one of its triangulations (b).
			The sides labeled as $e$ are to be identified, as well as the $0$-cells labeled as $v$.

			There are three kinds of points in the space $\DunceH$, and they are fully characterized by their own neighborhoods.
			If a point lies in the interior of the triangle, its neighborhood is certainly homeomorphic to $\mathbb{R}^2$.
			However, this is not going to happen to those points lying in the side $e$ nor to the special $0$-cell ``corner'' $v$.
			The boundary of a sufficiently small ball on $\DunceH$ centered at $v$ looks precisely like \Cref{fig:dunce_hat} (e) (the simplicial case is detailed in \Cref{stt:description_link_of_corner}).
			So, $v$ must be a vertex in any triangulation of $\DunceH$.
			Furthermore, any triangulation of $\DunceH$ subdivides $e$ into edges.

			For the forthcoming discussions, we consider \Cref{fig:dunce_hat} as a pictorial and notational reference.

			\subsubsection{Decomposition of $\DunceH$}
				\label{sec:decomposition_dunce_hat}

				To prove the last part of \Cref{stt:main_result}, we will follow the strategy of the previous sections, namely, we decompose $\DunceH$ into two partitionable spaces, then glue them back using our toolkit.
				Our decomposition will be as follows.

				Let $\TriangDH$ be a triangulation of $\DunceH$, with $v$ corresponding to the corner vertex.
				We use $v$ to split the facets of $\TriangDH$ into two subcomplexes:
				on the one hand, we have all the facets that include $v$, and on the other, those that do not.
				The corresponding subcomplexes are respectively the cone complex $v*\lk_{\TriangDH}(v)$ and the deletion complex $\del_{\TriangDH}(v)$ (see \Cref{fig:dunce_hat} (c)).
				Notice that the boundary shared by these two complexes is precisely $\lk_{\TriangDH}(v)$.

				We want to see what the link of $v$ looks like.
				The following well know result can be proved by an easy shelling argument.
				
				\begin{lemma}
					The link of a boundary vertex in any triangulation of the $2$-disk $\Disk$ (i.e.\ the ball $\Ball^2$) is a triangulation of $\Ball^1$ (i.e.\ a path).
					\label{stt:links_of_balls_as_balls}
				\end{lemma}

				\begin{remark}
					\Cref{stt:links_of_balls_as_balls} cannot be propagated to higher dimensions.
					For example, it is known that the double suspension of the Poincaré homology sphere is a $5$-sphere, but the link of a suspending vertex is not even a manifold (although it is a pseudomanifold with the same homologies of a sphere).
					\label{stt:on_poincare_homology_sphere}
				\end{remark}

				Recall that $e$ must be simplicially subdivided by edges.
				Hence, the corner vertex $v$ has to have distinct vertices $x$ and $y$ in its neighborhood lying on $e$ to ensure simpliciality.
				Thus, we can lift the triangulation $\TriangDH$ to a triangulation $\Delta$ of the triangle whose quotient space is $\DunceH$.
				To do that, we cut along $e$ as depicted in \Cref{fig:dunce_hat} (g).
				Observe that this action lifts copies of $v$, say $v_1, v_2$ and $v_3$, and also copies of $x$ and $y$.
				Denote as $\alpha, \beta$ and $\gamma$ to the respective quotients of $\lk_{\Delta}(v_1), \lk_{\Delta}(v_2)$ and $\lk_{\Delta}(v_3)$ once $\Delta$ is folded.

				\begin{lemma}
					Using the notation of the preceding paragraph, the subcomplexes $\alpha$ and $\beta$ triangulate $\Cycle$, while $\gamma$ triangulates $\Ball^1$.
					The vertex $x$ is common to $\alpha$ and $\gamma$, and the vertex $y$ to $\beta$ and $\gamma$.
					\label{stt:description_link_of_corner}
				\end{lemma}

				\begin{proof}
					\Cref{stt:links_of_balls_as_balls} ensures that $\lk_{\Delta}(v_i)$ is a path over the $1$-skeleton of $\Delta$.
					Given that the $x_i$'s and $y_i$'s are to be identified, respectively as $x$ and $y$, we get that $\alpha$ and $\beta$ become simplicial cycles and $\gamma$ the path joining them from $x$ to $y$.
					We observe that there is no common vertex between $\alpha$ and $\beta$, otherwise we would break simpliciality of $\TriangDH$.
				\end{proof}

				\begin{remark}
					Compare the link of $v$ described in \Cref{stt:description_link_of_corner} with the complex $L$ obtained in the first part of the proof of \cite[Theorem 4]{Zeeman:1964}; see also Figure 5 of the same paper.
				\end{remark}

				\begin{lemma}
					Let $\TriangDH$ be a triangulation of the dunce hat $\DunceH$, and let $v$ be its distinguished corner vertex.
					Then any triangulation of the complex $\left( \del_{\TriangDH}(v), \lk_{\TriangDH}(v) \right)$ is partitionable.
					\label{stt:twisted_underwear_partitionable}
				\end{lemma}

				\begin{proof}
					Since every edge in $\Delta$ belongs to at most two triangles, it follows that the complex $v_1 * \lk_{\Delta}(v_1)$ is a $2$-ball, and the link is a path (see \Cref{stt:links_of_balls_as_balls}).
					It follows by the Jordan-Schoenflies Theorem that the deletion of $v_1$ in $\Delta$ is a $2$-ball, since it is bounded by $\Cycle$.
					Use again the same argument for $v_2$ in the complex $\del_{\Delta}(v_1)$, and then for $v_3$ in the resulting complex.
					This procedure yields a triangulation of a disk relative to a connected part of its boundary.
					More precisely, in the way we labeled \Cref{fig:dunce_hat} (g), we end up with a triangulation of a disk relative to the path $\alpha' x_3 y_3 \gamma' x_2 y_2 \beta'$, where $\alpha', \gamma'$ and $\beta'$ are the link paths appearing after the iterative deletion of each $v_i$.
					A partitioning scheme of this disk is given by \Cref{stt:rel_disk_is_partitionable}.
					Then fold the complex with \Cref{stt:folding_lemma} to glue back the three copies of the path $xy$ lying on $e$.
					This yields a partitioning scheme of
					$\left( \del_{\TriangDH}(v), \lk_{\TriangDH}(v) \right)$.
				\end{proof}

				\begin{remark}
					Cutting $\TriangDH$ along $e$, as in the proof of \Cref{stt:twisted_underwear_partitionable}, proves that $\left( \TriangDH, \{v\} \right)$ is partitionable.
					Again, use the now-recurrent tandem of \Cref{stt:rel_disk_is_partitionable} and \Cref{stt:folding_lemma}.
					\label{stt:fun_facts_twisted_underwear_space}
				\end{remark}

			\medskip
			\phantomsection
				\label{sec:main_resutl_dunce_hat}

				Now we have the ingredients to finish the proof of \Cref{stt:main_result}.

				\begin{theorem}[\Cref{stt:main_result} for $\DunceH$]
					Any triangulation of $\DunceH$ is partitionable.
					\label{stt:dunce_hat_partitionable}
				\end{theorem}

				\begin{proof}
					We decompose a triangulation $\TriangDH$ of $\DunceH$ into $v*\lk_{\TriangDH}(v)$ and $\del_{\TriangDH}(v)$ as discussed earlier in this section.
					We know that
					$\left( \del_{\TriangDH}(v), \lk_{\TriangDH}(v) \right)$
					is partitionable by \Cref{stt:twisted_underwear_partitionable}.
					On the other hand, the cone $v*\lk_{\TriangDH}(v)$ is shellable hence partitionable:
					\Cref{stt:description_link_of_corner} tells us that $\lk_{\TriangDH}(v)$ is connected, hence shellable; therefore, the cone of $\lk_{\TriangDH}(v)$ is also shellable.

					Since
					$\lk_{\TriangDH}(v) = (v*\lk_{\TriangDH}(v)) \, \cap \, \del_{\TriangDH}(v)$,  \Cref{stt:assymetrical_gluing_lemma} gives us the desired partitionability of
					$\left(\left(v*\lk_{\TriangDH}(v)\right)
								\, \cup \,
								\del_{\TriangDH}(v),\
								\varnothing
					\right) = \TriangDH.$
					\qedhere
				\end{proof}


	\section*{Acknowledgments}
		\label{sec:acknowledgments}

		I am grateful to Russ Woodroofe for his guidance, ideas and careful reading of each draft of this paper.
		I also thank Bennet Goeckner, Masahiro Hachimori, Caroline Klivans, Jeremy Martin, Bruno Benedetti and Lorenzo Venturello for their helpful feedback and comments.
		I also thank the anonymous referees for their thoughtful comments. 
		
		The free open-source mathematics software system SageMath \cite{sagemath} and Hachimori’s online library of simplicial complexes \cite{Hachimori:web} were valuable resources.


	{
	\bibliography{../../../Master}
	\bibliographystyle{hamsplain}
	}

\end{document}